\documentclass[12pt]{article}

\usepackage{amssymb}

\makeatletter
\@addtoreset{equation}{section}
\makeatother

\newcommand{\beq}{\begin{equation}}
\newcommand{\eeq}{\end{equation}}
\newcommand{\bea}{\begin{eqnarray}}
\newcommand{\eea}{\end{eqnarray}}

\newcommand{\req}[1]{(\ref{#1})}
\newtheorem{thm}{Theorem}[section]

\newtheorem{cor}{Corollary}[section]
\newtheorem{prop}{Proposition}[section]

\begin{document}

\title{Cartesian approach for constrained mechanical system with three degree of freedom.}

\title{
{\bf Cartesian approach for constrained mechanical system with
three degree of freedom.
} \\
\vspace{1ex}
}

\date{}

\author{Rafael Ramirez,$^{1}$
Natalia Sadovskaia$^2$
\\
\vspace{0.001ex} \\
\normalsize  $^1$ Departament d'Enginyeria Inform\`{a}tica i Matem\`{a}tiques, \\
\normalsize  Universitat Rovira i Virgili, \\
\normalsize  Avinguda dels Pa\"{\i}sos Catalans 26, \\
\normalsize  43007 Tarragona, Spain \\
\normalsize  \\
\normalsize  $^2$Departamento de Matem\'atica Aplicada II, \\
\normalsize  Universitat Polit\`ecnica de Catalunya \\
\normalsize  C. Pau Gargallo,5 \\
\normalsize  08028, Spain \\
}

\maketitle

\begin{abstract}
In the history of mechanics, there have been two points of view
for studying mechanical systems: The Newtonian and the Cartesian.

According the Descartes point of view, the motion of mechanical
systems is described by the first-order differential equations in
the $N$ dimensional configuration space $\textsc{Q}.$

 In this paper we develop the Cartesian approach for mechanical
systems with three degrees of freedom and  with constraint which
are linear with respect to velocity. The obtained results we apply
to discuss the integrability of the geodesic flows on the surface
in the three dimensional Euclidian space  and to analyze the
integrability of a heavy rigid body in the Suslov and the Veselov
cases .
\end{abstract}

\section{Introduction.}
In "Philosophiae Naturalis Principia Mathematica" (1687), Newton
considers that movements of celestial bodies can be described by
differential equations of the second order. To determine their
trajectory, it is necessary to give the initial position and
velocity. To reduce the equations of motion to the investigation
of a dynamics system it is necessary to double the dimension of
the position space and to introduce the auxiliary phase space.

Descartes in 1644  proposed that the behavior of the celestial
bodies be studied from another point of view.  These ideas were
stated in "Principia Philosophiae" (1644) and in "Discours de la
m\'etode" (1637). According to Descarte the understanding of
cosmology starts from acceptance of the initial chaos, whose
moving elements are ordered according to certain fixed laws and
form the Cosmo. He consider that the Universe is filled with a
tenuous fluid matter (ether), which is constantly in a vortex
motion. This motion moves the largest particle of matter of the
vortex axis, and they subsequently form planets. Then, according
to what Descartes wrote in his "Treatise on Light", "the material
of the Heaven must be rotate the planets not only about the Sun
but also about their own centers...and this will hence form
several small Heavens rotating in the same direction as the great
Heaven."\cite{Koz1}. Thus the equation of motion in the Descartes
theory must be of the first order \begin{equation} \label{11}
  \dot{\textbf{x}}=\textbf{v}(x,t).\,
  \end{equation}
 Hence, to determine the trajectory from Descartes's point of
view it is necessary to give only the initial position.

In the modern scientific literature the study of the Descarte
ideas we can find in the monographic of V.V. Kozlov \cite{Koz1} in
which the author  give the following result.

\begin{thm}{\it The manifold $y=u(x,t),$ where $u$ is a covector on
$\textsc{Q}$ is an invariant manifold for the canonical
Hamiltonian equations with the Hamiltonian $H(x,y,t) $ if and only
if field $u$ satisfies the Lamb equation}
\begin{equation}\label{12}
    \partial_tu(x,t)+(rot u(x,t))v(x,t)=-grad h(x,t)
\end{equation}
 where $(rotu)=\partial_xu-(\partial_xu)^t$ ia skew-symmetric
$n\times{n}$ matrix,
\begin{equation}\label{13}
   v(x,t)=\partial_yH(x,y)|_{y=u(x,t)},\quad
h(x,t)=H(x,y,t)|_{y=u(x,t)}
\end{equation}
\end{thm}

From the physical standpoint, equations \req{11},\,\req{12} and
\req{13}describe the motion of the collisionless medium: particles
moving along different trajectories do not interact.

In \cite{Koz1} affirm that  "solving dynamics problem is possible
inside the configuration space".  For this it is necessary to
solve Lamb equations which is a system of partial differential
equations on $\textsc{Q}$, and then, using \req{13} to calculate
the vector field $v$ from the solution of the Lamb equation to
solve \req{11}.

In \cite{Ram1} we developed the Cartesian approach for mechanical
system with configuration space $\textsc{Q}$ and with constraints
linear with respects to velocity. The aim of the present paper is
to develop the results obtained  for mechanical system with three
degrees of freedom in the particular case in which
$\textsc{Q}=\mathbb{E}^3$ is the three dimensional Euclidean space
and $\textsc{Q}=S0(3)$ is the special orthogonal group of
rotations of  $\mathbb{E}^3.$


\section{Cartesian vector
 field on three dimensional Euclidean space }

Let $\mathbb{E}^3$ be the three dimensional Euclidian space with
cartesian coordinates $x=(x_1,x_2,x_3).$

We consider a particle with Lagrangian function

$$
L=\frac{1}{2}||\dot{\textbf{x}}||^2-U(x)$$ and constraints
 \beq (\dot{\textbf{x}},\textbf{a})=0\label{21}\eeq where $(\,,\,)$ denotes the
 scalar product in $\mathbb{E}^3,$
$\dot{{\textbf{x}}}=(\dot{x}_1,\,\dot{x}_2,\,\dot{x}_3)$ and
${\textbf{a}}(x)=(a_1(x)\,,a_2(x),\,a_3(x))$ is a smooth vector
field in
 $\mathbb{E}^3$

It is well known that the equations of motions can be deduced from
the d'Alembert-Lagrange principle \cite{Ram2}
\begin{equation}\label{22}
\left\{%
\begin{array}{ll}
    \ddot{{\textbf{x}}}=U_\textbf{x}+\mu {\textbf{a}}(x),  \\
    (\textbf{a},\dot{\textbf{x}})=0,, \\
\end{array}%
\right.
\end{equation}
where $\mu$ is the Lagrangian multiplier,
$U_\textbf{x}=(U_{x_1},\,U_{x_2},\,U_{x_3}),\,
U_{x_j}=\partial_{x_j}U.$

In \cite{Ram1} we introduce the following definition

{\textbf{ Definition 1}}

{\it We say the smooth vector field $\textbf{v}
(x)=(v_1(x),\,v_2(x),\,v_3(x))$ is the Cartesian vector field for
a constrained particle in $\mathbb{E}^3$ with the constraints
\req{21} if}
 \begin{equation}
 \label{23}
 [\textbf{v}(x)\times{rot{\textbf{v}(x)}}]=\Lambda (x)\textbf{a}(x).\,\end{equation}
 where
 $[\,\times\,]$  denotes the vector product in
$\mathbb{E}^3,$
$$rot\textbf{v}=(\partial_{x_2}v_3-\partial_{x_3}v_2,\,\partial_{x_3}v_1-\partial_ {x_1}v_3,\,\partial_{x_1}v_2-\partial_ {x_2}v_1).$$
and $\Lambda$ is a function:
$$\Lambda=\displaystyle\frac{1}{||\textbf{a}||^2}([{\textbf{v}(x)}\times rot{{\textbf{v}(x)}}],{\textbf{a}(x)}).$$

By a simple computation from \req{23} we can see that

\begin{equation}
\label{24}
\left\{%
\begin{array}{ll}
    (\textbf{a}(x),\textbf{v}(x))=0,  \\
    (\textbf{a}(x),rot{\textbf{v}(x)})=0.\\
\end{array}%
\right.   \end{equation}

\begin{cor}

{\it Let $\textbf{v}$ the Cartesian vector field. Then the
following relations hold}
\begin{equation}
\label{25}
\left\{%
\begin{array}{ll}
\ddot{{\textbf{x}}}=(\frac{1}{2}||\textbf{v}(x)||^2)_x+\Lambda (x){\textbf{a}}(x)\\
(\dot{\textbf{x}},\textbf{a})=0,\end{array}%
\right.   \end{equation}
\end{cor}
The proof it is easy to obtain in view of the equality
$$\ddot{\textbf{x}}=(\frac{1}{2}||\textbf{v}(x)||^2)_x+ [\textbf{v}(x)\times{rot{\textbf{v}(x)}}]$$
which is deduced after derivation the differential equations
generated by the vector field $\textbf{v}$

\begin{equation}
\label{26}
  \dot{\textbf{x}}=\textbf{v}(x)\end{equation}
 The  system \req{25} can be obtained from the Lagrangian
 equations with
Lagrangian function
$$\textsc{L}=\frac{1}{2}||\dot{\textbf{x}}-\textbf{v}(x)||^2$$
where $\textbf{v}$ is a Cartesian vector field.

{\textbf{Definition 2}}

{\it The study of the behavior of the constrained particle in
$\Bbb{E}^3$ by using the equations \req{22} or \req{26},\req{23}
or \req{25} say the Classical, Cartesian and Lagrangian approach
respectively.}

We illustrate the above ideas in the following example

{\textbf{A non-holonomically constrained particle in
$\textsc{R}^3.$}}

Consider a particle with the kinetic energy
$T=\frac{1}{2}||\dot{\textbf{x}}||^2$ and non-holonomic
constraints
$$\dot{x}_1+\hat{a}(x_3)\dot{x}_2=0$$

 This instructive academic example, in the particular case when $\hat{a}(x_3)=x_3$ due to Rosenberg
\cite{Ros}. This example was also used to illustrate the theory in
Bates and Sniatycki \cite{Bates}.

 The Descartes approach in this case
produces the vector field $\textbf{v}:$
$$ {\textbf{v}}=\lambda_2(\hat{a}(x_3)\partial_{x_1}-\partial_{x_2}) -
{\lambda}_3\partial_{x_3}$$ and condition \req{24}
 for this case takes the form
$$(rot{{{\textbf{v}}}},\textbf{a}(x))=0\Longleftrightarrow\frac{1}{2}\partial_{x_3}((1+\hat{a}^2){\lambda}^2_2)+
(\hat{a}(x_3)\partial_{x_1}{\lambda}_3-\partial_{x_2}{\lambda}_3){\lambda}_2=0.$$

 We shall study the case when this relation holds in view of the
 equalities
$${\lambda}_2=\displaystyle\frac{A}{\sqrt{1+\hat{a}^2(x_3)}},\quad {\lambda}_3=b_2(x_3),$$ for
$A$ an arbitrary constant and $b_2$ an arbitrary function on
$x_3.$

The equations generated by the vector field ${\bold{v}}$ in this
case can be written as
\begin{equation}\label{A}
  \left\{%
\begin{array}{ll}
   \dot{x_1}=\displaystyle\frac{\hat{a}(x_3)A}{\sqrt{1+\hat{a}^2(x_3)}} \\
\dot{x_2}=-\displaystyle\frac{A}{\sqrt{1+\hat{a}^2(x_3)}}\\
\dot{x_3}=-b_2(x_3)
 \end{array}%
\right.
\end{equation}

 The all trajectories of
these equations are easy to obtain.

The Lagrangian approach produces the following differential
equations

\[
  \left\{%
\begin{array}{ll}
\ddot{x}_1=-b(x_3)\partial_{x_3}(\displaystyle\frac{A\hat{a}(x_3)}{\sqrt{1+\hat{a}^2(x_3)}})\\
\ddot{x}_2=b(x_3)\partial_{x_3}(\displaystyle\frac{A}{\sqrt{1+\hat{a}^2(x_3)}})\\
\ddot{x}_3=\partial_{x_3} \displaystyle\frac{1}{2}b^2(x_3)\end{array}%
\right.
\]

 \begin{cor}

 All the trajectories of the equation of motion of the
constrained Lagrangian system
$$<\, \textsc{E}^3,\,L=\frac{1}{2}||\dot{\textbf{x}}||^2-
U(x_3),\,\{\dot{x}_1+\hat{a}(x_3)\dot{x}_2=0\}\,>$$ can be
obtained from \req{A} with $b(x_3)=\pm{\sqrt{h+U(x_3)}}.$
\end{cor}

In fact, the equations of motion obtained from the
D'Alembert-Lagrange Principle are
\[
  \left\{%
\begin{array}{ll}
\ddot{x}_1=\mu\\
\ddot{x}_2=\hat{a}(x_3)\mu\\
\ddot{x}_3=\partial_{x_3}U(x_3)\\
\dot{x}_1+\hat{a}(x_3)\dot{x}_2=0 \end{array}%
\right.
\]

Therefore,
$$\frac{d}{dt}(\dot{x}_2-\hat{a}(x_3)\dot{x}_1)=-\frac{d\hat{a}(x_3)}{dx_3}\dot{x}_3\dot{x}_1$$
hence,
$$\dot{x}_2=\displaystyle\frac{A}{\sqrt{1+\hat{a}^2(x_3)}}$$
where $A$ is an arbitrary constant.

On the other hand from the
equation$$\ddot{x}_3=\partial_{x_3}U(x_3)$$ we easily obtain
$\dot{x}_3=\mp\displaystyle\sqrt{2(U(x_3)+h)},$ where $h$ is an
arbitrary constant.

Finally by considering the constraints we deduce the system of the
first order ordinary differential equations \req{A}. In this
example the Descartes, the lagrangian and Classical approach
coincide .

Below we determine the Cartesian vector field for a particle on
the surface in  $\mathbb{E}^3.$

 First we introduce the ve ctor fields  $X,\,,Y,\,Z$ which are characteristic elements of the
1-form
$$\Omega =a_1(x)dx_1+a_2(x)dx_2+a_3(x)dx_3$$
\begin{equation}
\label{27}
\left\{%
\begin{array}{ll}
X=a_3\partial_y-a_2\partial_z\nonumber\\
Y=a_1\partial_z-a_3\partial_x\\
Z=a_2\partial_x-a_1\partial_y\end{array}%
\right.   \end{equation}

Clearly,   the most general element of the given 1-form $\Omega$
is

\[ \textbf{v}=w_1X+w_2Y+w_3Z\]
hence

\begin{equation}
\label{28}  {\textbf{v}}(x)=[{\textbf{a}}(x)\times
{\textbf{w}}(x)]
\end{equation}

  where $\textbf{w}(x)=(w_1(x),\,w_2(x),\,w_3(x))$ is
an arbitrary smooth vector field which we shall determine in such
a way that \req{24} takes place.

By using the identity
$$ rot[\textbf{a}(x)\times{\textbf{b}(x)}]=[\textbf{a},\,\textbf{b}]+div{\textbf{b}(x)}\textbf{a}(x)-div{\textbf{a}(x)}\textbf{b}(x)$$
where $[\textbf{a},\,\textbf{b}]$ is the Lie bracket of the smooth
vector field $\textbf{a}$ and $\textbf{b},$ one can prove the
following assertion

\begin{cor}
The condition \req{23},\req{28}  can be written as follows
\begin{eqnarray}
 \label{29}
 div([\textbf{a}(x)\times [\textbf{a}(x)\times{\textbf{w}(x)}]])=([\textbf{a}(x)\times rot \textbf{a}(x)],\textbf{w}(x))\end{eqnarray}
\end{cor}
\begin{prop}

{\it Let us suppose that the vector field $\textbf{a}$ is such
that
\begin{eqnarray}
 \label{30}[\textbf{a}(x)\times rot\textbf{a}(x)]=\textbf{0}\end{eqnarray}
 then the Cartesian vector field exist i and only if
\begin{eqnarray}
 \label{31}
\textbf{a}(x)=f_\textbf{x}(x)\end{eqnarray} for a certain smooth
function $f.$}
\end{prop}
{\textbf{Proof}}

From \req{29},\,\req{30} follows that
$$[\textbf{a}(x)\times
[\textbf{a}(x)\times{\textbf{w}(x)}]]$$ is a solenoidal vector
field, hence
$$[\textbf{a}(x)\times [\textbf{a}(x)\times\textbf{w}(x)]]=rot\textbf{W}(x)$$
for arbitrary vector field ${\textbf{W}} ,$ thus  the following
representation holds
$$\textbf{w}=\frac{(\textbf{a}(x),\textbf{w}(x))}{||\textbf{a}(x)||^2}\textbf{a}(x)+\frac{rot\textbf{W} (x)}
{||\textbf{a}(x)||^2}$$ as a consequence
\begin{eqnarray}
 \label{32}
 (\textbf{a}(x),{rot}{\textbf{W}} (x) )=0 .\end{eqnarray}

Clearly if $({\textbf{a}} (x) ,{rot}{\textbf{W}} (x) )\ne{0}$ then
the Cartesian vector field does not exist if \req{30} holds. If we
choose
$$\textbf{a}(x)=f_\textbf{x}(x),\quad \textbf{W} (x)=\Phi G_\textbf{x}(f,\Phi )$$
when  $\Phi,\,G$ are an arbitrary smooth functions, then  we
obtain that \req{32} holds identically and a consequence the
vector $\textbf{w}$ takes the form
\begin{eqnarray}
 \label{33}
\textbf{w}(x)=\frac{(f_\textbf{x}(x),\textbf{w}(x))}{||f_\textbf{x}(x)||^2}+\nu
(x)[f_\textbf{x}(x)\times{\Phi_\textbf{x}}],\quad \nu
(x)=\frac{\partial_fG(f,\Phi
)}{||f_\textbf{x}(x)||^2}\end{eqnarray}

\begin{cor}

{\it The Cartesian vector field for a particle in $\mathbb{E}^3$
which is constrained to move on the surface
\begin{eqnarray}
 \label{34}
 f(x)=c,\,c\ne{0}\end{eqnarray}
 generated the following differential system}
\begin{eqnarray}
 \label{35}
 \dot{\textbf{x}}=\nu (x)(||f_\textbf{x}(x)||^2\Phi_\textbf{x}(x)-(f_\textbf{x}(x),\,\Phi_\textbf{x}(x))f_\textbf{x}(x))\end{eqnarray}

\end{cor}
\begin{cor}

{\it The Lagrangian approach for a particle in $\mathbb{E}^3$
which is constrained to move on the surface \req{34} produces the
following differential equations}
\begin{eqnarray}
 \label{36}
 \ddot{\textbf{x}}=\frac{\partial}{{\partial \textbf{x}}}\Big(\frac{1}{2}\nu^2 (x)||f_\textbf{x}(x)||^2||[f_\textbf{x}(x)\times
 \Phi_\textbf{x}(x)]||^2\Big)+\Lambda (x)f_\textbf{x}(x)\end{eqnarray}
\end{cor}
\begin{cor}

{\it If there exist a function $G$ and $\Phi$ such that
\begin{equation}
\label{37} ||[f_\textbf{x}(x)\times
 \Phi_\textbf{x}(x)]||^2=\displaystyle\frac{2h(f)g}{G_f(f,\Phi )}\equiv{\Psi (f,\Phi )\,g}
\end{equation}

Then the equations \req{36} take the form}
\end{cor}
\begin{eqnarray}
 \label{38}
 \ddot{\textbf{x}}=\lambda_0(x)f_\textbf{x}(x),\quad \lambda_0(x)=h_f(f)+\Lambda (x).\end{eqnarray}
If one introduce the matrix $A(x):$
\begin{equation}\label{39}
A(x)=\left(
\begin{array}{lll}
f_{x_1x_1}&f_{x_1x_2}&f_{x_1x_3}\\
f_{x_1x_2}&f_{x_2x_2}&f_{x_2x_3}\\
f_{x_1x_3}&f_{x_2x_3}&f_{x_3x_3},\end{array}\right)\end{equation}then
one checks, that the equations \req{36} may be written as
$$\ddot{x}=\
\frac{(A(x)\textbf{v}(x),\textbf{v}(x))}{||f_\textbf{x}(x)||^2}f_\textbf{x}(x),$$
where $\textbf{v}$ is the Cartesian vector field generated the
differential equations \req{35}.  The differential equations
\req{38} determined the geodesic flows on the surface \req{34} and
admits the energy integral
$$||\dot{\textbf{x}}||^2=2h(f).$$
If there is an additional first integral, functionally independent
with the energy integral , then the geodesic flow is integrable.

In order to study the integrability of the geodesic flow on the
given surface we introduce the following functions which we
determine from \req{35}
\[
\label{39}
\left\{%
\begin{array}{ll}
    F_1=\displaystyle{\Big(\frac{||[f_\textbf{x}\times\dot{\textbf{x}}]||}{||f_\textbf{x}||}\Big)^2}\Longleftrightarrow{F_1=||\dot{x}||^2},
     \\
     \\
    F_2=\displaystyle{\Big(\frac{||f_\textbf{x}||||[\Phi_\textbf{x}\times\dot{{\textbf{x}}}]||}{(f_\textbf{x},\,\Phi_\textbf{x})}\Big)^2},
     & \hbox{if $(f_\textbf{x},\,\Phi_\textbf{x})\ne{0} $;} \\
     \\
     F_3=\displaystyle{\Big(\frac{||\Phi_\textbf{x}||||[\textbf{x}\times\dot{{\textbf{x}}}]||}{||[\textbf{x}\times\Phi_\textbf{x}]||}\Big)^2},
     & \hbox{if $(f_\textbf{x},\,\Phi_\textbf{x})=0,\quad\Phi_\textbf{x}\ne{\kappa (x) \textbf{x}}.$} \\
\end{array}%
\right.
\]

In view of \req{37}, it is easy to show that
$$F_j=2h(f),\quad j=1,2,3.$$

\section{Integrability of the geodesic flow on the homogeneous surface. }

We now consider the surface
\begin{equation}
 \left\{%
\begin{array}{ll}
    f(x)=c,\quad c\ne{0}, \\
    (\textbf{x},f_\textbf{x}(x))=mf(x). \\
\end{array}%
\right.
\end{equation} which we will call the
{\it homogeneous surface of degree $m.$}

From the Euler formula follow that $c=0$ is the unique critical
value of $f.$ hence for $c\ne{0}$ the function
$$g=||f_\textbf{x}(x)||^2>0$$
on the given homogeneous surface.

 Taking (3.1) into account we deduce the relations
\begin{eqnarray}
\label{42}
 A(x)\textbf{x}^{T}&=&(m-1)f^{T}_\textbf{x}(x),\\
(\textbf{x},g_\textbf{x}(x))&=&2(m-1)g(x)\end{eqnarray}

Below we use the following notation \[\{F,G,H\}\equiv\left|
\begin{array}{rrr}
F_{x_1}&F_{x_2}&F_{x_3}\\
G_{x_1}&G_{x_2}&G_{x_3}\\
H_{x_1}&H_{x_2}&H_{x_3}\end{array}\right|\]

Clearly, if $F,\,G,\,H$ are independent functions then
$\{F,G,H\}\ne{0}.$

 The integrability of
the geodesic flow on the homogeneous surface we shall study in the
following two cases
\begin{eqnarray}
\label{43}
\{f,g,r^2\}&=&0, \\
\{f,g,r^2\}&\ne&0
\end{eqnarray}
where $r^2=x^2_1+x^2_2+x^2_3.$

 We analyze the first case. We study only the particular subcase when
the homogeneous surface is such that
\begin{eqnarray}
\label{44} ||f_\textbf{x}(x)||^2=g(f,r).\end{eqnarray} Hence, in
view of (3.3) we give
\begin{eqnarray}
\label{45} mf\partial_fg+r\partial_rg=2(m-1)g(f,r).\end{eqnarray}
 We assume that the arbitrary function $\Phi$ is such that
 $$\Phi_\textbf{x}=\textbf{x},$$
thus the differential equation generated by the Cartesian vector
field and second order differential equations of the geodesic
flows under the indicated condition take the form respectively
\begin{eqnarray}
\label{46} \dot{\textbf{x}}=\nu (x)(g\,\textbf{x}-m\,ff_\textbf{x}),\\
\ddot{\textbf{x}}=\frac{m\partial_rg f
h(f)}{r^2g^2}f_\textbf{x}\end{eqnarray} where

\begin{equation}\label{46_1}
 \nu^2g(f,r)(g(f,r)r^2-m^2f^2)
\end{equation}

 \begin{prop}

{\it The geodesic flow on the homogeneous surface under the
assumption \req{44} is integrable}
\end{prop}
{\textbf{Proof}}

First we observe that there is the function $\nu$ such that
\req{46_1} holds , i.e.,
\[
G^2_f(f,r )=\displaystyle\frac{2h(f) g(f,r)}{g(f,r)r^2-m^2f^2},\]

 hence  exist the additional first
integral $F_2$ which in this case takes the form
$$F_2=\frac{g(f,r)||[f_\textbf{x}(x)\times\dot{{\textbf{x}}}]||^2}{m^2f^2}\Leftrightarrow g(f,r)
||[f_\textbf{x}(x)\times\dot{\textbf{x}}]||^2=2{m^2f^2}h(f).$$

The particular class of the study homogeneous surface are the
following.

If $m=1$ then $(\textbf{x},g_\textbf{x})=0,$ in particular this
relation holds if
$$g=\Psi (\frac{f}{r})$$
A concrete example we obtain from the celestial mechanics
\cite{Dub}:
\begin{eqnarray}
\label{47}
 f(x)=r+(\textbf{b},\textbf{x})=c,\quad c\ne{0},\end{eqnarray}

  where $\textbf{b}=(b_1,b_2,b_3)$ is
a constant vector field. In this case we have
$$g=\frac{2f}{r}+||b||^2-1$$
The first integral for this particular case are

\begin{equation}\left\{%
\begin{array}{ll}
    ||\dot{\textbf{x}}||^2=2h(f), \\
    (\displaystyle\frac{2f}{r}+||\textbf{b}||^2-1)||[\textbf{x}\times\dot{{\textbf{x}}}]||^2=2f^2h(f) \\
\end{array}%
\right.
\end{equation}

It is interesting to deduce the  equations of motion of a particle
constrained to move in the surface \req{47} with the subsidiary
condition that there is a nonzero constant vector field
$\breve{c}=(c_1,c_2,c_3):$
\begin{eqnarray}
\label{48}
 (\textbf{x},\breve{\textbf{c}})=0,\quad (\textbf{b},\check{\textbf{c}})=0\quad \Rightarrow (f_\textbf{x}(x),\breve{\textbf{c}})=0,\end{eqnarray}
from the Lagrangian approach.

 By choosing the function $\Phi$ as follows
$$\Phi=(c_3-c_2)x_1+(c_1-c_3)x_2+(c_2-c_1)x_3$$
and introducing the new time $\sigma$ as
$$d\sigma=(\displaystyle\frac{x_1+x_2+x_3}{r}+b_1+b_2+b_3)dt$$
and letting the prime denote differentiation with respect to
$\sigma,$ we have that the equation generated by the Cartesian
vector field can be written as
\begin{eqnarray}
\label{49}
{\textbf{x}}^{'}=\frac{1}{||\breve{\textbf{c}}||^2}[f_\textbf{x}(x)\times{\breve{\textbf{c}}}]\end{eqnarray}
By considering that

$$rot[f_\textbf{x}\times{\breve{\textbf{c}}}]=\frac{\breve{\textbf{c}}}{r}$$ we obtain that the Lagrangian approach generated
the second order differential equations
$${\textbf{x}}^{''}=-\frac{f(\textbf{x})x}{||\breve{\textbf{c}}||^2 r^3}$$
Thus , if $$f(x)=||\breve{\textbf{c}}||^2$$ then we obtain the
well known equations
\begin{eqnarray}
\label{50}{\textbf{x}}^{''}=-\frac{\textbf{x}}{ r^3}\end{eqnarray}

These equations admit the following first integrals

\begin{equation}
\left\{%
\begin{array}{ll}
  ||{\textbf{x}}^{'}||^2 = \displaystyle\frac{2}{r}+\frac{||\textbf{b}||^2-1}{||\breve{\textbf{c}}||^2} \\
  \\
 \left[{\textbf{x}}^{'}\times \breve{\textbf{c}} \right]= -(\displaystyle\frac{\textbf{x}}{r}+\textbf{b}) \\
\\
 \left[\textbf{x}\times{{\textbf{x}}}^{'}\right] = \breve{\textbf{c}}\\
\end{array}%
\right.
\end{equation}

These relations are easy to obtain from \req{50}.

The  equations \req{50}, after the orthogonal transformation
$$\xi=\frac{(\textbf{b},\textbf{x})}{||\textbf{b}||},\quad \eta=\frac{([
\breve{\textbf{c}}\times{\textbf{b}}],\textbf{x})}{||\breve{\textbf{c}}||\,||\textbf{b}||},\quad
\zeta=\frac{(\breve{\textbf{c}},\textbf{x})}{||\breve{\textbf{c}}||},$$take
the form
\begin{equation}
\left\{%
\begin{array}{ll}
    {\xi}^{''}=-\displaystyle\frac{\xi}{\sqrt{(\xi^2+\eta^2)^3}},\\
    \\
    {\eta}^{''}=-\displaystyle\frac{\eta}{\sqrt{ (\xi^2+\eta^2)^3}} \\
    \zeta=0. \\
\end{array}%
\right.
\end{equation}
These equations describe the behavior of the particle with
Lagrangian function
$$L=\frac{1}{2}({\xi^{'}}^2+{\eta^{'}}^2)-\frac{1}{\sqrt{ (\xi^2+\eta^2)}}$$

 constrained to move on the one curve of the family of conics
$$\breve{f}(x)=\sqrt{\xi^2+\eta^2}+||b||\xi=||\breve{c}||^2$$
The differential equations generated by the Cartesian vector field
in this coordinates can be represented in Hamiltonian form with
Hamiltonian function $\breve{f}$ \cite{Ram3}

 We now turn to the study the particular case of the homogeneous
 surface with the condition $g=g(f,r).$

 If  $$g=r^{2(m-1)}\Psi (\frac{f}{r^m}),\quad \Psi(\frac{f}{r^m})\ne{(\frac{f}{r^m})^2} $$ then after the change
 $$F=\int\frac{d(\xi )}{\sqrt{\Psi (\xi )-\xi^2}},\quad
 \xi=\frac{f}{r^m},$$
 we obtain
$$||F_\textbf{x}(x)||^2=\frac{1}{r^2}.$$

Finally, if
$$g=f^{\frac{2(m-1)}{{m}}}\Psi (\frac{f}{r^m})$$
then after the change
$$f=F^m$$
we deduce the equation
$$||F_\textbf{x}(x)||^2=\tilde{\Psi}(\frac{F}{r}),$$
which show that this case is equivalent to the first case study
above.

We have already studied the case in which $\{f,g,r^2\}=0.$ Now we
begin to study the case in which the functions $f,\,g,\,r^2$ are
independent. Hence

\begin{equation}\label{50}
  \{f,g,r^2\}\ne{0}.
\end{equation}
Under this assumption we obtain that
$$x_j=x_j(f,g,r^2),\quad j=1,2,3$$
thus we deduced that
\begin{equation}\label{51}
\left\{%
\begin{array}{ll}
\Phi=\Phi (f,g,r^2)\\
\Phi_\textbf{x}=\partial_f\Phi\,f_\textbf{x}+\partial_g\Phi\,g_\textbf{x}+\displaystyle\frac{\partial_r\Phi}{r}\,\textbf{x}\\
\end{array}%
\right.
\end{equation}


\begin{prop}

{\it If there exists the functions $\Phi$ and $G$ such that}
\begin{equation}\label{52}
\left\{%
\begin{array}{lll}
\Phi=\Phi (f,g,r^2),\quad G=G(f,\Phi )\\
g\nu=G_f(f,\Phi )\\
 ||f_\textbf{x}||^2||\Phi_\textbf{x}\times f_\textbf{x}||^2\nu^2=2h(f)
\end{array}%
\right.
\end{equation}
{\it then the geodesic flow on the homogeneous surface of degree
$m>1$ is integrable.}
\end{prop}
{\textbf{{Proof}}

We prove this assertion only for the case when

\begin{equation}\label{B}
(\Phi_\textbf{x},\,f_\textbf{x})=0,
\end{equation} thus the surface
$\Phi=c_1$ is orthogonal to the given homogeneous surface.  Under
this assumption we obtain that the differential equations deduced
from the Cartesian and Lagrangian approach are respectively
\begin{equation}\label{53}
\left\{%
\begin{array}{ll}
\dot{\textbf{x}}=\nu (x)\Phi_\textbf{x}(x)\\
\ddot{\textbf{x}}=\lambda_0(x)f_\textbf{x}(x)\end{array}%
\right.
\end{equation} where $\lambda_0$ can
be determined as follows
$$\lambda_0(x)=\nu^2(f_\textbf{x},\partial_\textbf{x}(\frac{1}{2}||\Phi_x||^2))=h_f(f)-
\frac{\nu||\Phi_\textbf{x}||^2}{g}(f_\textbf{x},\nu_\textbf{x})$$
After derivation the function
$$F_3=\Big(\frac{||\Phi_\textbf{x}||||[\textbf{x}\times\dot{{\textbf{x}}}]||}{||[\textbf{x}\times\Phi_\textbf{x}]||}\Big)^2$$
along the solutions of \req{53} we obtain
$$\displaystyle\frac{d{F}_3}{dt}=(\Phi_\textbf{x},\partial_\textbf{x}(\nu^2 ||\Phi_\textbf{x}||^2)).$$
which is equal to zero  in view of \req{52},\,\req{B} thus the
function $F_3$ is the first integral of the geodesic flow.

Clearly, in order to assess the integrability of the geodesic flow
in this case we need first to check whether function $\nu$ exists
such that \req{52} holds.

\section{The geodesic flow on the quadrics and the third-order surface }

In order to illustrate the above ideas we consider the algebraic
surface of degree three:

\begin{equation}\label{54}
  f(x)=x_1\,x_2\,x_3=c,\quad c\ne{0}.
\end{equation}

This case was examined already by Riemann in his study of motion
of a homogeneous liquid ellipsoid. More exactly, Riemann examined
the integrability of the geodesic flow on \req{54}.

In \cite{Koz2} the author state the following problem.

"Is it true that the geodesic flow on a generic third-order
algebraic surface  is not integrable?. In particular I do not know
a rigorous proof of non-integrability for the surface \req{54}"

By considering that in this case
$$g=(x_1x_2)^2+(x_1x_3)^2+(x_3x_2)^2$$
thus the functions $f,g$ and $r^2$ are independent. The dependence
$x_j=x_j(f,g,r^2),\,j=1,2,3$ we obtain as follows.

We introduce the cubic polynomial in $z:$
$$P(z)=z^3-r^2z^2+g z-f^2=(z-x^2_1)(z-x^2_2)(z-x^2_3),$$
and by using  Cardano's formula we obtain the require dependence.

In order to construct the Cartesian approach in this case first we
observe that the surface
$$\Phi (\xi , \eta ,\zeta )=c_1$$
where
$$\xi=\frac{1}{2}(x^2_1-x^2_2),\quad\eta=\frac{1}{2}(x^2_3-x^2_1),\quad
\zeta=\frac{1}{2}(x^2_2-x^2_3)$$ is orthogonal to surface
\req{54}. Thus the differential equations generated by the
cartesian vector field are
\begin{equation}\label{55}
\left\{%
\begin{array}{lll}
\dot{x}_1=\nu(\Phi_\xi-\Phi_\eta )\,x_1\\
\dot{x}_2=\nu(\Phi_\zeta-\Phi_\xi )\,x_2\\
\dot{x}_3=\nu(\Phi_\eta-\Phi_\zeta )\,x_3
\end{array}%
\right.
\end{equation}
To determine the existence the solution of \req{52} or, what is
the same, \[
||\Phi_\textbf{x}||^2=\displaystyle\frac{2h(f)}{G^2_f(f,\Phi
)}\equiv{\Psi (f,\Phi )},\]
  is for us an open problem.

Now we study the subcase when the given surface is such that

\begin{equation}\label{56}
  f(x)=\frac{1}{2}(b_1x^2_1+b_2x^2_2+b_3x^2_3)
\end{equation}

First we state and solve the following problem.

{\textbf{Problem 1}}

Let $X,\,Y,\,Z$ are the vector fields \req{27},\,\req{31}.

We require to determine the function $f$ in such a way that these
vector field formed a three dimensional Lie algebra.

The solution of this problem it is easy to obtain in view of the
equality
\begin{equation}\label{57}\Upsilon_1=A(x)\,\Upsilon_2, \end{equation}
where $A$ is the matrix given by the formula \req{39} and
$$\Upsilon_1=col([Y,Z],\,[Z,X],\,[X,Y]),\quad
\Upsilon_2=col(X,\,Y,\,Z)$$ and by using the Bianchi
representation
\begin{equation}\label{58}\Upsilon_3=B(x)\,\Upsilon_4, \end{equation}
where $\Upsilon_3=col([U,V],\,[V,W],\,[W,U]),\quad
\Upsilon_4=col(U,\,V,\,W),$ where $U,\,V,\,W$ are the vector
fields, $B$ is the matrix:

\[B(x)=\left(
\begin{array}{lll}
0&a&b_3\\
b_1&0&0\\
0&b_2&-a\end{array}\right)\] and $a,\,b_1,\,b_2,\,b_3$ are certain
constants

 From \req{57} and \req{58} after integration we obtain the class
 of functions which generated the three dimensional Lie algebra.

\begin{equation}\label{59}
\left\{%
\begin{array}{llllll}
 1.\quad f=b_1\,x^2+b_2\,y^2+b_3\,z^2\\
2.\quad f=b_1\,x^2+a\,(y^2-z^2)+2\,b\,yz\\
3.\quad  f=2b\,y\,x+b_3\,z^2\\
4.\quad  f=b\,y^2+2\,b_1\,z\,x\\
\end{array}%
\right.
\end{equation}

We construct the Cartesian vector field for the first case.

In view of the relation
$$g=b^2_1x^2_1+b^2_2x^2_2+b^2_3x^2_3$$
we observe that $f,\,g,\,r^2$ are independent functions. The
following equalities it is easy to obtain:

\begin{equation}
\left\{%
\begin{array}{llll}
    x^2_1=\displaystyle{\frac{b_2b_3r^2-2(b_2+b_3)f+g}{(b_1-b_2)(b_1-b_3)}}\\
& &\\
     x^2_2=\displaystyle{\frac{b_1b_3r^2-2(b_1+b_3)f+g}{(b_2-b_1)(b_2-b_3)}} \\
& &\\
     x^2_3=\displaystyle{\frac{b_2b_1r^2-2(b_2+b_1)f+g}{(b_3-b_2)(b_3-b_1)}}\\
\end{array}%
\right.\end{equation}

  Notice that
$$\Phi (\xi ,\eta , \zeta )=c_1,$$  where
$$\xi=\frac{x^{b_2}_3}{x^{b_3}_2},\quad
\eta=\frac{x^{b_1}_2}{x^{b_2}_1},\quad
\zeta=\frac{x^{b_3}_1}{x^{b_1}_3}$$is an orthogonal surface to the
given surface we obtain that the differential equations generated
by the Cartesian vector field in this case can be written as
\begin{equation}\label{0}
\left\{%
\begin{array}{llll}
   \dot{x}_1=\displaystyle{\frac{\nu}{x_1} (\Phi_\zeta \zeta b_3-\Phi_\eta \eta b_2)}\\
&&\\
   \dot{x}_2=\displaystyle{\frac{\nu}{x_2} (\Phi_\eta \eta b_1-\Phi_\xi \xi b_3)} \\
&&\\
    \dot{x}_1=\displaystyle{\frac{\nu}{x_3} (\Phi_\xi\xi b_2-\Phi_\zeta\zeta b_1)} \\
\end{array}%
\right.\end{equation}

Now we introduce the elliptic coordinates in $\mathbb{R}^3:$
\begin{equation}
\left\{%
\begin{array}{ll}
x^2_1=\displaystyle\frac{(\lambda_1+b^{-1}_1)(\lambda_2+b^{-1}_1)(\lambda_3+b^{-1}_1)}{(b^{-1}_1-b^{-1}_2)(b^{-1}_1-b^{-1}_3)}\\
x^2_2=\displaystyle\frac{(\lambda_1+b^{-1}_2)(\lambda_2+b^{-1}_2)(\lambda_3+b^{-1}_2)}{(b^{-1}_2-b^{-1}_1)(b^{-1}_2-b^{-1}_3)}\\
x^2_3=\displaystyle\frac{(\lambda_1+b^{-1}_3)(\lambda_2+b^{-1}_3)(\lambda_3+b^{-1}_3)}{(b^{-1}_3-b^{-1}_2)(b^{-1}_3-b^{-1}_1)}\\
\end{array}%
\right.\end{equation} where $\lambda_1,\,\lambda_2,\,\lambda_3$
are the roots of the cubic polynomial in $w:$
\[
\left\{%
\begin{array}{llllll}
-w^3+(x_3^2+x_1^2+x_2^2-b^{-1}_1-b^{-1}_2-b^{-1}_3)w^2+(x_1^2(b^{-1}_2+
b^{-1}_3)+x_2^2(b^{-1}_1+b^{-1}_3)\\
+x_3^2(b^{-1}_1+b^{-1}_2)-b^{-1}_1b^{-1}_2-b^{-1}_1b^{-1}_3-b^{-1}_2b^{-1}_3)w\\
+b^{-1}_1b^{-1}_2b^{-1}_3(b_1x_1^2+b_2x_2^2+b_3x_3^2-1)=0\end{array}%
\right.\]

In this coordinates we obtain
$$||\dot{x}||^2=g_{11}\dot{\lambda}^2_1+g_{33}\dot{\lambda}^2_2+g_{33}\dot{\lambda}^2_3$$
where
\begin{equation}
\left\{%
\begin{array}{ll}
g_{11}=\displaystyle\frac{(\lambda_1-\lambda_2)(\lambda_1-\lambda_3)}{4(\lambda_1-b^{-1}_1)(\lambda_1-b^{-1}_2)\lambda_1-b^{-1}_3)}\\
g_{22}=\displaystyle\frac{(\lambda_2-\lambda_1)(\lambda_2-\lambda_3)}{4(\lambda_2-b^{-1}_1)(\lambda_2-b^{-1}_2)\lambda_2-b^{-1}_3)}\\
g_{33}=\displaystyle\frac{(\lambda_3-\lambda_2)(\lambda_3-\lambda_3)}{4(\lambda_3-b^{-1}_1)(\lambda_3-b^{-1}_2)\lambda_3-b^{-1}_3)}\\
\end{array}%
\right.\end{equation} The differential equations \req{0} in
elliptic coordinates can be transformed to the form

\[
\left\{%
\begin{array}{llll}
   \dot{w}=0\\
&&\\
   \dot{u}=\displaystyle\frac{2\nu}{b_1b_2b_3}((b_2-b_3)\Phi_\zeta\zeta+
   (b_1-b_2)\Phi_\eta \eta+(b_3-b_1)\Phi_\xi \xi)\equiv\Psi_1 (\lambda_1,\lambda_2) \\
&&\\
    \dot{v}= \displaystyle\frac{2\nu}{b_1b_2b_3} (\displaystyle\frac{1}{b_1}+\displaystyle\frac{1}{b_2}+\displaystyle\frac{1}{b_3})({(b_2-b_3)}b_1\Phi_\zeta\zeta+(b_1-b_2)b_3\Phi_\eta \eta+
   (b_3-b_1)b_2\Phi_\xi \xi)\equiv\Psi_2 (\lambda_1,\lambda_2)\\
u=\lambda_1+\lambda_2+\lambda_3,\quad
v=\lambda_1\lambda_2+\lambda_1\lambda_3+\lambda_2\lambda_3,\quad
w=\lambda_1\lambda_2\lambda_3,
\end{array}%
\right.\]and by putting $$\lambda_3=0,$$  after some calculations
we deduce the planar system
\[
\left\{%
\begin{array}{llll}
   \dot{\lambda}_1=\displaystyle\frac{\Psi_1 (\lambda_1,\lambda_2)\lambda_1-\Psi_2 (\lambda_1,\lambda_2))}
   {\lambda_1-\lambda_2}\\
\dot{\lambda}_2=\displaystyle\frac{\Psi_1
(\lambda_1,\lambda_2)\lambda_2-\Psi_2 (\lambda_1,\lambda_2))}
   {\lambda_1-\lambda_2}.
\end{array}%
\right.\]
 In order to deduce the differential equations
for the  geodesic flow by using the Lagrangian approach first, it
is necessary in the first place obtain the solution of the
equations \req{52}.

The integrability of the geodesic flow on the quadric (m=2) by
using the classical approach, was proved by Jacobi and Chasles.

\section{The geometrical and physical meaning of the Cartesian
vector field}

The purpose of this section is to determine the geometrical and
physical meaning of the Cartesian vector field constructed above.

 Hertz's principle of least
curvature is a special case of Gauss' principle, restricted by the
two conditions that there be no applied forces and that all masses
are identical. (Without loss of generality, the masses may be set
equal to one.) Under these conditions, Gauss' minimized quantity
can be written
\[\textsc{Z}=\sum_{j=1}^N |\displaystyle\frac{d^2x^j}{dt^2}|^2\]

The kinetic energy \[T=\frac{1}{2}||\dot{\textbf{x}}||^2\] is also
conserved under these conditions

Since the line element $ds^2$ in the $3N$-dimensional space of the
coordinates is defined

\[ds^2=2Tdt^2\Longleftrightarrow\displaystyle\frac{ds^2}{dt^2}=2T \]
by considering the conservation of energy we obtain
\[\displaystyle\frac{ds^2}{dt^2}=2h\]

Dividing $\textsc{Z}$ by $2T$ yields another minimal quantity

\[\textsc{K}=\sum_{j=1}^N |\displaystyle\frac{d^2x^j}{ds^2}|^2\]
Since $\sqrt{\textsc{K}}$ is the local curvature of the trajectory
in the $3N$-dimensional space of the coordinates, minimization of
$K$ is equivalent to finding the trajectory of least curvature (a
geodesic) that is consistent with the constraints. Hertz's
principle is also a special case of Jacobi's formulation of the
least-action principle. Curvature refers to a number of loosely
related concepts in different areas of geometry.  In mathematics,
a geodesic is a generalization of the notion of a straight line to
curved spaces. Definition of geodesic depends on the type of
curved space. If the space carries a natural metric then geodesics
are defined to be (locally) the shortest path between points on
the space.

Below we restricted to the case when the configuration space is
the three dimensional Euclidean space with Cartesian coordinates
$x=(x_1,\,x_2,\,x_3).$  The geodesic flow on the surface $f(x)=c$
is determined by the second- order differential equations
\[\displaystyle\frac{d^2x^j}{dt^2}=\mu f_{x_j},\quad \mu =
\displaystyle\frac{(A(x)
\dot{\textbf{x}},\,\dot{\textbf{x}})}{||f_\textbf{x}||^2} \quad
j=1,2,3\]

which, by considering the energy integral, can be written as
follows
\[\displaystyle\frac{d^2x^j}{ds^2}=\tilde{\mu} f_{x_j},\quad\tilde{\mu} =\displaystyle\frac{(A(x) \tau,\,\tau )}{||f_\textbf{x}||^2} \quad j=1,2,3\]

where
$${\tau}=\displaystyle\frac{d\textbf{x}}{ds},\quad ||{\tau}||^2=1.$$
Clearly that
$$\sqrt{\textsc{K}}=\displaystyle\frac{|(A(x) {\tau},\,{\tau} )|}{||f_\textbf{x}||^2}$$

 The Hertz's Principle of Least Curvature and problem on the
 determination the principal directions on the surface lead us to state the following problem.

\textbf{Problem 2}

Determine the
$$extremum (A(x){\tau},{\tau})$$
under the conditions
\[
\left\{%
\begin{array}{ll}
    ||\tau||^2-1=0 \\
    (f_\textbf{x},\tau)=0\\
\end{array}%
\right.\]

\textbf{Solution}

Note that in this case the Lagrangian function is

$$\textsc{L}=(A(x)\tau,\tau)+\sigma (f_\textbf{x},\tau)+z(||\tau||^2-1)$$
where $\sigma$ and $z$ are the Lagrangian multiplier and computer
\[\left\{%
\begin{array}{lllll}
    \displaystyle{\frac{\partial\textsc{L}}{\partial\tau_j}}=0,\quad j=1,2,3
    \Longleftrightarrow{ (A(x)+zI)\tau^T+\sigma f^T_\textbf{x}=0}\\
    &&\\
  \displaystyle{\frac{\partial\textsc{L}}{\partial\sigma}}=0\Leftrightarrow{(f_\textbf{x},\tau)=0} \\
&&\\
    \displaystyle{\frac{\partial\textsc{L}}{\partial z}}=0\Leftrightarrow{ ||\tau||^2-1=0}\\
\end{array}%
\right.\] where $\tau^T=col(\tau_1,\,\tau_2,\,\tau_3)$ and $I$ is
the diagonal matrix: $I=diag{(1,\,1,\,1)},$ from the first group
of equations, we deduced the following equalities
\begin{equation}\label{60}
    \tau^T(x)=-\sigma (A(x)+zI)^{-1}f^T_\textbf{x}(x))
\end{equation}
and
\begin{equation}\label{B1}
 \textsc{R}_z\chi=\vec{0}
\end{equation}
if $$\det{(A(x)+zI)}\ne{0},$$ where
$\chi=col(\tau_1,\,\tau_2,\,\tau_3,\,\sigma )$ and $\textsc{R}_z$
is the following family of matrixes
\[ \textsc{R}_z=\left(%
\begin{array}{cccc}
  f_{x_1x_ 1}+z& f_{x_1x_2 } & f_{x_1x_3 } & f_{x_1} \\
  f_{x_1x_2 } & f_{x_2x_2 }+z & f_{x_2x_3 } & f_{x _2} \\
  f_{x_1x_ 3} & f_{x_2x_ 3} & f_{x_3x_ 3}+z& f_{x_3} \\
  f_{x_1} & f_{x_2} & f_{x_3}  & 0 \\
\end{array}%
\right)\]

In view of that the vector $\chi$ is non-zero vector then from
\req{B1} one can deduce  that

\begin{equation}\label{62}
 \det{\textsc{R}_z}=0
\end{equation}

In order to establish the relation between the vector field with
components given by \req{60} and Cartesian vector field \req{35}
we  introduce the family of vector fields $\textbf{v}_z:$

\begin{equation}\label{62}
  \textbf{v}_z=((A(x)+zI)^{-1}f_x,\partial_x),\quad
  \partial_x=(\partial_{x_1},\,\partial_{x_2},\,\partial_{x_3})
\end{equation}

where $z$ is a complex parameter. After some calculations one can
prove that $v_z$ admits the  representations

\begin{equation}\label{63} \textbf{v}_z=\left|
\begin{array}{cccc}
  f_{x_1x_ 1}+z& f_{x_1x_2 } & f_{x_1x_3 } & f_{x_1} \\
  f_{x_1x_2 } & f_{x_2x_2 }+z & f_{x_2x_3 } & f_{x _2} \\
  f_{x_1x_ 3} & f_{x_2x_ 3} & f_{x_3x_ 3}+z& f_{x_3} \\
  \partial_{x_1} & \partial_{x_2} & \partial_{x_3}  & 0 \\
\end{array}%
\right|=\displaystyle{\frac{z^2X_1+zX_2+X_3}{\det{(A+zI)}}}\end{equation}
which are equivalents equivalent to \req{62},
 where $\triangle{f}=\partial_{x_1x_1}f+\partial_{x_2x_2}f+\partial_{x_3x_3}f.$
and $X_1,\,X_2,\,X_3$ denote the vector fields:
\begin{equation}\label{64}
\left\{%
\begin{array}{lllll}
  X_1=(f_x,\partial_x) \\
 X_2=(g_x,\partial_x)-\triangle{f}\,(f_x,\partial_x),\quad g=||f_x||^2 \\
  X_3=\det{A(x)}\,v_{z}|_{z=0}
\end{array}%
\right.
\end{equation}

We now introduce the function
\begin{equation}\label{65}
  F(z)=df(v_z)=\displaystyle{\frac{z^2df(X_1)+zdf(X_2)+df(X_3)}{\det{(A+zI)}}}
\end{equation}
which in view of \req{63} may be written as
$$F(z)=\det{(A+zI)}\det{\textsc{R}_z}$$

\begin{cor}
{\it If the function $f$ is a homogeneous function of degree $m$
then}

\begin{equation}\label{66}
X_3=g^2(x)K(x)\,(x,\,\partial_x),\quad g=||f_x||^2
\end{equation}
\end{cor}
where $K$ is the Gaussian curvature of the homogeneous surface
which one can calculate as follows
\begin{equation}
\label{67}
K(x)=\left\{%
\begin{array}{ll}
    \displaystyle{ -\frac{f\{f_{x_1},\,f_{x_2},\,f_{x_3}\}}{(m-1)g^2}} & \hbox{if $m\ne{1} $;} \\
   \displaystyle{-\frac{f\{f,\,f_{x_2},\,f_{x_3}\}}{x_1\,g^2}} & \hbox{if $m={1} $;} \\
\end{array}%
\right.\end{equation}

 Below we shall study only the case when the function $f:$

\begin{equation}\label{68}
(df(X_2))^2-4df(X_3)df(X_1)>0
\end{equation}
 Clearly, under this assumption the function $F$ has two different real
 roots which we shall denote by $z_1,\,z_2$

\begin{cor} {\it Let $\textbf{v}_1,\,\textbf{v}_2$ are the vector fields
such that
$$\textbf{v}_j=\textbf{v}_z|_{z=z_j},\quad j=1,2 $$

then the solution of the problem 2 are the vector fields :}
\end{cor}
\[\left\{
\begin{array}{c}
 \tau^{(j)}=\displaystyle{\frac{\textbf{v}_j}{||\textbf{v}_j(x)||}},\quad j=1,2 \\
\\
  (\tau^{(1)},\,\tau^{(2)})=0
\end{array}\right.\]
The proof it is easy to obtain.

\begin{prop}

{\it Let $f(x)=c,\,c\ne{0}$ be the homogeneous surface which
satisfies \req{68}.

Then the most general vector field tangent to the given surface
admits the development}
\begin{equation}\label{69}
 \textbf{v}(x)=[f_\textbf{x}(x)\times{[f_\textbf{x}(x)\times{(\mu_1(x)\,g_\textbf{x}(x)+\mu_2(x)\,\textbf{x}\,
 )}}]].
\end{equation}
where $\mu_1,\,\mu_2 $ are arbitrary smooth functions.
\end{prop}
\textbf{Proof}

One can check direct from the above that the most general vector
field tangent to the given homogeneous surface can be written as
$$\textbf{v}(x)=a_1(x)\textbf{v}^{(1)}(x)+a_2(x)v^{(2)}(x)$$
where $a_1,\,a_2$ are arbitrary smooth functions.

A brief calculation show that

\begin{equation}\label{70}
 v(x)=\lambda_1(x)\,f_\textbf{x}(x)+\lambda_2(x)\,g_\textbf{x}(x)+\lambda_3(x)\,\textbf{x}
\end{equation}
where

\[\left\{%
\begin{array}{lll}
   \lambda_1(x)=\hat{a}_1(x)(z^2_1-\triangle{f}\,z_1)+\hat{a}_2(x)(z^2_2-\triangle{f}\,z_2) \\
   \lambda_2(x)=\hat{a}_1(x)\,z_1+\hat{a}_2(x)\,z_2 \\
    \lambda_3(x)=g^2(x)K(\hat{a}_1(x)+\hat{a}_2(x))\\
 a_j(x)=\hat{ a}_j\det{(A(x)+z_jI)},\quad j=1,2
\end{array}%
\right.\]

One can see that the equation

\begin{equation}\label{71}
 g\lambda_1+(f_\textbf{x},\,g_\textbf{x})\lambda_2+mf\lambda_3=0
\end{equation}
holds identically.

To complete the proof, we show the equivalence of \req{69} and
\req{70}, \req{71}. Indeed, using \req{71} we obtain that
$$\lambda_1=-(\frac{(f_\textbf{x},\,g_\textbf{x})\lambda_2+mf\lambda_3}{g}),$$
inserting into \req{70} and introducing the notations
$$\mu_1(x)=\frac{\lambda_3}{g},\quad
\mu_2(x)=\frac{\lambda_2}{g}$$ we get \req{69}.

\begin{prop}

{\it The vector field \req{69} is Cartesian vector field if the
following relation holds}
\begin{equation}\label{72}
  \{f,g,r\}(\partial_r\mu_1-\partial_g\mu_2-\displaystyle{\frac{\mu_2}{g}})=0
\end{equation}
\end{prop}
From the definition we obtain that the given vector field is
Cartesian (see definition 1) if the following equality takes place

\begin{equation}
\label{73} (
f_\textbf{x},(rot([f_\textbf{x}(x)\times{[f_\textbf{x}(x)\times{(\mu_1(x)\,g_\textbf{x}(x)+\mu_2(x)\,\textbf{x}\,)
 )}}]])=0
\end{equation}

which is equivalent to \req{72}.

From \req{72} after straightforward calculations we can prove the
following assertion

\begin{cor}
{\it Let us suppose that \req{72} holds, then the vector field
\req{70} admits the representation}
\end{cor}

\begin{equation}\label{74}
 \textbf{v}(x)=\kappa \,[f_\textbf{x}\times [f_\textbf{x}\times \Phi_\textbf{x}]]
\end{equation}
 where $\kappa $ and
$\Phi $ functions such that
\begin{equation}
\kappa= \label{75}
\left\{%
\begin{array}{ll}
\kappa (f,\Phi) & \hbox{if $\{f,g,r\}\ne{0} $;} \\
\\
\mu_1\,\displaystyle{\frac{g_r}{r}}+\mu_2& \hbox{if $\{f,g,r\}={0} $;} \\
\end{array}%
\right.\end{equation} and
\begin{equation}
\Phi_x= \label{76}
\left\{%
\begin{array}{ll}
\Phi_\textbf{x}(f,g,r) & \hbox{if $\{f,g,r\}\ne{0} $;} \\
\\
\textbf{x} & \hbox{if $\{f,g,r\}={0} $;} \\
\end{array}%
\right.\end{equation}

By comparing \req{74}, \req{75} with \req{35} we  obtain that the
solution of the problem 2 under the condition \req{75},\,\req{76}
coincide with the Cartesian vector field constructed in the
section 3. In such a way  we obtain the physical and geometrical
meaning of the constructed above Cartesian vector field.

\textbf{Remark.}

The physical and geometrical meaning of the Cartesian vector field
for the case when the given vector field $\textbf{a}:$

\begin{equation}\label{L}
(\textbf{a},rot\textbf{{a}})\ne{0}\end{equation} can be obtained
analogously to the case study above   by considering that under
condition \req{L}  the equation deduced from the Lagrangian
approach, can be written as follows
$$\ddot{\textbf{x}}=\displaystyle\frac{(A(x)
\dot{\textbf{x}},\,\dot{\textbf{x}})}{||f_\textbf{x}||^2}\textbf{a}$$
where
\[A(x)=\left(
\begin{array}{lll}
\partial_1a_1&\frac{1}{2}(\partial_1a_2+\partial _2a_1) &\frac{1}{2}(\partial_1a_3+\partial _3a_1)\\
\frac{1}{2}(\partial_1a_2+\partial _2a_1)&\partial_2a_2&\frac{1}{2}(\partial_2a_3+\partial _3a_2)\\
\frac{1}{2}(\partial_1a_3+\partial
_3a_1)&\frac{1}{2}(\partial_2a_3+\partial
_3a_2)&\partial_3a_3\end{array}\right)\]

The problem 2 in this case we can state as follows

\textbf{Problem 3}

Determine the
$$extremum (A(x){\tau},{\tau})$$
under the conditions
\[
\left\{%
\begin{array}{ll}
    ||\tau||^2-1=0 \\
    (\textbf{a},\tau)=0\\
\end{array}%
\right.\]

\section{Descartes approach for non-holonomic system
 with three degree of freedom  and
one constraints .}

 Our goal in this section is to extend the
Cartesian approach developed above for natural mechanical system
with configuration space $$\textsc{Q},\quad \dim\textsc{Q}=3$$ in
this space the metric (kinetic energy)
$$T=\frac{1}{2}\sum_{j,k=1}^3G_{jk}(x)\dot{x}^j\dot{x}^k\equiv{\frac{1}{2}||\dot{x}||^2}$$
allows calculating the $rot$ of the vector field $v$ on
$\textsc{Q}.$ The invariant definition of $rot{v}$ we can find in
\cite{Koz1}.  If we assign a covector field $p=(p_1,\,p_2,\,p_3)$
with components

\begin{equation}\label{77}
  p_j=\sum_{k=1}^3G_{jk}(x)v^k{(x)}
\end{equation}
to the vector field $v=(v^1(x),\,v^2(x),\,v^3(x)),$ then the
components of $rot{v}$ we can write explicitly
$$rot{v(x)}=(\frac{1}{\sqrt{G}}(\partial_2p_3-\partial_3p_2),\,
\frac{1}{\sqrt{G}}(\partial_3p_1-\partial_1p_3),\frac{1}{\sqrt{G}}(\partial_1p_2-\partial_2p_1))$$
where $G=\det{(G_{kj}(x))}.$

The vector field \req{28} for this mechanical system we shall
represented as follows \begin{equation}\label{78}
   \dot{\textbf{x}}=[\textbf{a}(x)\times{(\lambda_1\,\textbf{b}(x)+\lambda_2\textbf{c}(x))}]
\end{equation}
where $\textbf{a},\,\textbf{b},\,\textbf{c}$ are the independent
smooth vector in $\textsc{Q},$ i.e.,
\begin{equation}\label{79}\Upsilon= (\textbf{a},[\textbf{b}\times
\textbf{c}])\ne{0}\end{equation} and $\lambda_1,\,\lambda_2$ are
smooth function which we determine as a solution of the equation
\begin{equation}\label{80}
(\textbf{a},rot{[\textbf{a}\times{(\lambda_1\,\textbf{b}+\lambda_2\textbf{c})}]})=0\end{equation}

The Lagrangian  approach produces the following second-order
differential equations
\begin{equation}\label{81}
\displaystyle\frac{d}{dt}\displaystyle\frac{\partial
T}{\partial\dot{x}^k}-\displaystyle\frac{\partial
T}{\partial{x}^k}=\displaystyle\frac{\partial
\frac{1}{2}||\textbf{v}||^2}{\partial{x}^k}+\Lambda a_k(x),\quad
k=1,2,3
\end{equation}

 We shall illustrate this case for the Chapliguin-Caratheodory sleigh and for the heavy rigid body in the
 Suslov case.

\textbf{The Chapliguin-Caratheodory sleigh }

We shall now analyze one of the most classical nonholonomic
systems : Chapliguin-Carathodory's sleigh \cite{NF}. The idealized
sleigh is a body that has three points of contact with the plane.
Two of them slide freely but the third, $A,$ behaves like a knife
edge subjected to a constraining force $\mathbb{R}$ which does not
allow transversal velocity.  More precisely, let $yoz$ be an
inertial frame and $\xi\,A\eta$ a frame moving with the sleigh.
Take as generalized coordinates the Descartes coordinates of the
center of mass $C$ of the sleigh and the angle $x$ between the $y$
and the $\xi$ axis. The reaction force $\mathbb{R}$  against the
runners is exerted laterally at the point of application $A$ in
such a way that the $\eta$ component of the
 velocity is zero. Hence, one has the constrained system
 $\textsc{{M}}$ with the configuration space
$X=S^1\times\mathbb{R}^2,$ with the kinetic energy
$$T=\frac{m}{2}(\dot{y}^2+\dot{z}^2)+\frac{I_c}{2}\dot{x}^2,$$ and
with the constraint
$$\epsilon\dot{x}+\sin{x}\dot{y}-\cos{x}\dot{z}=0,$$ where $m$ is
the mass of the system and $J_c$ is the moment of inertia about a
vertical axis through $C$ and $\epsilon=|AC|.$ Observe that the
"javelin" (or arrow or Chapliguin's skate) is a particular case of
this mechanical system and can be obtained when $\epsilon=0$

To apply the Descartes approach for this system, first we
determine the vector $\textbf{b}$ and $\textbf{a}$ in such a way
that the determinant $\Upsilon\ne{0}.$
 In this subcase, we achieve this condition if
$$\textbf{a}=(\epsilon,\,\sin x,\,-\cos x)\quad
\textbf{b}=(0, \cos{x},\,sin x ),\quad
       \textbf{c}=(1,0,0).$$
Under these restrictions we obtain that $\Upsilon=1$ and it is
easy to show that the vector field $\textbf{{v}}$ takes the form:
$$\textbf{v}=\lambda_3(\partial_x+\epsilon\sin x\partial_y+\epsilon \cos x\partial_z)-\lambda_2 (\cos x\partial_y-
\sin x\partial_z). $$

The Descartes  approach  produce the differential equations
\cite{Sad}

\begin{equation}\label{K0}
\left\{%
\begin{array}{lll}
\dot{x}=\lambda_3(x,y,z,\epsilon )\\
\dot{y}=\lambda_2(x,y,z,\epsilon )\cos{x}-
\epsilon \lambda_3\sin{x}\\
\dot{z}=\lambda_2(x,y,z,\epsilon )\sin{x}+\epsilon\lambda_3\cos{x}
\end{array}%
\right.\end{equation}
 where $\lambda_2,\,\lambda_3$ are solutions
of the partial differential equations
\begin{equation}\label{K1}
  \sin x(J\partial_z\lambda_3+\epsilon
m\partial_y\lambda_2)+\cos x(J\partial_y\lambda_3-\epsilon
m\partial_z\lambda_2)-m(\partial_x\lambda_2-\epsilon\lambda_3)=0
\end{equation}
where $J=J_C+\epsilon^2m.$

Clearly,
$$||\textbf{v}||^2=(J_C+m\epsilon^2)\lambda^2_3(x,y,z,\epsilon\,)+m\lambda^2_2(x,y,z,\epsilon\,)$$

 Hence, for the arrow ($\epsilon=0$) we
have
\[\left\{%
\begin{array}{ll}
\dot{x}=\lambda_3(x,y,z,0 )\\
\dot{y}=\lambda_2(x,y,z,0 )\cos{x}\\
\dot{z}=\lambda_2(x,y,z,0 )\sin{x}
\end{array}%
\right.\]
$$J_C(\sin x\partial_z\lambda_3+\cos
x\partial_y\lambda_3)-m\partial_x\lambda_2=0$$

Clearly, the equation \req{K1} holds in particular if
\[\left\{%
\begin{array}{ll}
\partial_y\lambda_3=\displaystyle\frac{\epsilon m}{J_C+\epsilon^2m}\partial_z\lambda_2\\
\partial_z\lambda_3=-\displaystyle\frac{\epsilon m}{J_C+\epsilon^2m}\partial_y\lambda_2\\
\partial_x\lambda_2=\epsilon\lambda_3
\end{array}%
\right.\]

 After some calculations we can prove that the functions

\[\left\{%
\begin{array}{ll}
 \lambda_2=\cos\alpha\,V_1(y,z,\epsilon )+
\sin{\alpha}\,V_2(y,z,\epsilon )+
a\int{K(x,\epsilon )dx}\\
\lambda_3=\displaystyle\frac{am}{J_C+a^2m}\Big(\/cos{\alpha}\,V_2(y,z,a)-
\sin{\alpha}\,V_1(y,z,\epsilon )\Big)+K(x,\epsilon ),\\
\alpha=\displaystyle\frac{\epsilon^2mx}{J_C+\epsilon^2m}
\end{array}%
\right.\]

are solutions of \req{K1}, where  $K$ is an arbitrary function and
$V_1,\,V_2$ are functions which satisfy the  Cauchy-Riemann
conditions:
\[\left\{%
\begin{array}{ll}
\partial_yV_1(y,z,\epsilon )&=\partial_zV_2(y,z,\epsilon )\\
\partial_zV_1(y,z,\epsilon )&=-\partial_yV_2(y,z,\epsilon ).\end{array}%
\right.\]

\begin{cor}
{\it The all trajectories of the Chapliguin skate ($\epsilon=0$)
under the action of the potential field of force with components
$(0,\,mg,\,0)$ can be obtained from the Descartes approach}
\end{cor}
\textbf{Proof.}

In fact, for the case when $\epsilon=0$ the classical approach for
Chapliguin-Carathodory's sleigh gives the following equations of
motion
\[\left\{%
\begin{array}{ll}
\ddot{x}=0\\
\ddot{y}=mg+\sin x\mu \\
\ddot{z}=-\cos x\mu\\
\sin x\dot{y}-\cos x\dot{z}=0\end{array}%
\right.\]
 Hence,by derivation  we obtain
$$\frac{d}{dt}(\frac{\dot{z}}{\sin x})=g\cos x$$
as a consequence,
\[\left\{%
\begin{array}{ll}
\dot{x}=C_0,\quad C_0\ne{0}\\
\dot{y}=(\displaystyle\frac{g\sin x}{C_0}+C_1)\cos x\\
\dot{z}=(\displaystyle\frac{g\sin x}{C_0}+C_1)\sin x
\end{array}%
\right.\]
 or,

\[\left\{%
\begin{array}{ll}
\dot{x}=0\\
\dot{y}=({gt\cos x_0}+C_1)\cos x_0\\
\dot{z}=(gt\cos x_0+C_1)\sin x_0
\end{array}%
\right.\]
 Clearly, the solutions of these
equations coincide with the solutions of \req{K0},\,\req{K1}  with
the subsidiary conditions \cite{Sad}
$$J_C\lambda^2_3+m\lambda^2_2=mgy+h.$$

 are particular cases of the equations obtained from the
Descartes approach.
\begin{cor}
{\it The all trajectories of  Chapliguin -Carathodory's sleigh by
inertia can be obtained from the Descartes approach}
\end{cor}
\textbf{Proof}

Let us suppose that
$$\lambda_j=\lambda_j(x,\epsilon),\quad j=1,2$$
then the all trajectories of the equation \req{K0} can be obtained
from the formula
\[\left\{%
\begin{array}{ll}
y=y_0+\int\displaystyle\frac{(\lambda_2(x,\epsilon )\cos{x}-
\epsilon \lambda_3\sin{x})dx}{\lambda_3(x,\epsilon )}\\
{z}=z_0-\int\displaystyle\frac{(\lambda_2(x,y,z,\epsilon
)\sin{x}-\epsilon\lambda_3\cos{x})dx}{\lambda_3(x,\epsilon )}\\
t=t_0+\int\displaystyle\frac{dx}{\lambda_3(x,\epsilon )}
\end{array}%
\right.\]

 On the other hand, for the Chapliguin- Caratheodory
sleigh by inertia from the classical approach we deduce the
following equations
\[\left\{%
\begin{array}{ll}
J_C\ddot{x}=\epsilon\mu\\
m\ddot{y}=\sin x\mu \\
m\ddot{z}=-\cos x\mu\\
\epsilon\dot{x}+\sin x\dot{y}-\cos x\dot{z}=0
\end{array}%
\right.\]

 Hence, after
straightforward calculations we obtain the system
\[\left\{%
\begin{array}{ll}
\dot{x}=qC_0\cos (q\epsilon x+C),\quad q^2=\displaystyle\frac{m}{J_C+m\epsilon^2}\\
\dot{y}=C_0(\sin (q\epsilon x+C)\cos x-q\epsilon\cos (q\epsilon x+C)\sin x)\\
\dot{z}=C_0(\sin (q\epsilon x+C)\sin x+q\epsilon\cos
(q\epsilon+C)\cos x)
\end{array}%
\right.\] where $ q^2=\displaystyle\frac{m}{J_C+m\epsilon^2},$

 which are particular case of the equations
\req{K0} with
$$\lambda_2=C_0\sin (q\epsilon x+C),\quad \lambda_3=C_0q\cos (q\epsilon x+C)$$
Evident that in this case
$$2||\bold{v}||^2=(J_C+m\epsilon^2)\lambda^2_3(x,\epsilon\,)+m\lambda^2_2(x,\epsilon\,)\equiv{mC^2_0}$$

\section{ The rigid body around a fixed point in the Suslov and
Veselov cases.}

 In this section we study one classical problem of
non-holonomic dynamics formulated by Suslov \cite{Koz3}. We
consider the rotational motion of a rigid body around a fixed
point and subject to the non-holonomic constraints
$(\tilde\textbf{a},\omega)=0$ where ${\omega}$ is a body angular
velocity and $\tilde\textbf{a}$ is a constant vector. Suppose the
body rotates in an force field with potential
$U(\gamma_1,\gamma_2,\gamma_3)$. Applying the method of Lagrange
multipliers we write the equations of motion in the form
\begin{equation}
\label{81}
\left\{%
\begin{array}{ll}
  I\dot{{\omega}}= [I\omega\times{\omega}]+[{\gamma}\times\frac{\partial
U}{\partial{\gamma}}]+\mu \tilde\textbf{a} \\
 \dot{{\gamma}}=[{\gamma}\times{\omega}] \\
  (\tilde\textbf{a},\bold{\omega})=0
\end{array}%
\right.
\end{equation}
Where
\[
\begin{array}{c}
  I=diag(I_1,I_2,I_3) \\
  {\gamma}=(\gamma_1=\sin z\sin x,\quad\gamma_2=\sin z\cos
x,\quad\gamma_3=\cos z)
\end{array}
\]
 $I_1,\,I_2,\,I_3$ are the
inertial moment of the body.

 If we assume that the vector
$\tilde\textbf{a}=(0,0,1)$ \cite{Koz3}, then

\begin{equation}
\label{82}
\left\{%
\begin{array}{ll}
I_1\dot{\omega}_1=\gamma_3\partial_{\gamma_2}U-\gamma_2\partial_{\gamma_3}U\\
I_2\dot{\omega}_2=\gamma_1\partial_{\gamma_3}U-\gamma_3\partial_{\gamma_1}U\\
(I_1-I_2)\omega_1\omega_2+\gamma_2\partial_{\gamma_1}U-\gamma_1\partial_{\gamma_2}U+\mu=0\\
\dot{\gamma}_1=-\gamma_3\omega_2\\
\dot{\gamma}_2=\gamma_3\omega_1\\
\dot{\gamma}_3=\gamma_1\omega_2-\gamma_2\omega_1
\end{array}%
\right.
\end{equation}
The above system always has two independent first integrals
\[
\begin{array}{ll}
K_1=\frac{1}{2}(I_1\omega^2_1+I_2\omega^2_2)-U(\gamma_1,\gamma_2,\gamma_3)\\
K_2=\gamma^2_1+\gamma^2_2+\gamma^2_3
\end{array}%
\]
 For the real
motions $K_2=1.$

 By the Jacobi's theorem about the last
multiplier, if there exits a third independent first integral
$K_3$ which is functionally independent together with $K_1$ and
$K_2,$ then the Suslov problem is integrable  by quadratures
\cite{Koz3}

To determine the integrable cases of the Suslov problem seems
interesting the following result which we can prove after
straightforward calculations.

\begin{prop}

{\it Let us suppose that the  potential function $U$ in \req{82}
is determine as follows

\begin{equation}\label{83}
  U=\frac{1}{2I_1I_2}(I_1\mu^2_1+I_2\mu^2_2)-h
\end{equation}
 where $\mu_1,\,\mu_2$are solutions of the partial
differential equations

\begin{equation}\label{84}
\gamma_3(\frac{\partial \mu_1}{\partial\gamma_2}- \frac{\partial
\mu_2}{\partial \gamma_1})-\gamma_2\frac{\partial
\mu_1}{\partial\gamma_3}+\gamma_1\frac{\partial \mu_2}{\partial
\gamma_3}=0,
\end{equation}
 then the equations \req{82}, \req{83}
  admits the first integrals }
\begin{equation}\label{85}
 I_1\omega_1=\mu_2,\quad
I_2\omega_2=-\mu_1
\end{equation}
\end{prop}
 The aim of this apartat is to propose the Descartes approach for heavy rigid body in the Suslov
case.

Let us suppose that $\textsc{Q}=SO(3),$ with the Riemann metric
\[G=\left(
\begin{array}{ccc}
  I_3 & I_3\cos{z} & 0 \\
  I_3\cos{z} & (I_1\sin^2 x+I_2\cos^2 x)\sin^2{z}+I_3\cos^2{z} & (I_1-I_2)\sin x\cos x\sin z \\
  0 & (I_1-I_2)\sin x\cos x\sin z & I_1\cos^2x+I_2\sin^2x\\
\end{array}\right)
\]

$$\det G=I_1I_2I_3\sin^2z,$$

In this case we have that the constraints are
$$\omega_3=0\Leftrightarrow\dot{x}+\cos z \dot{y}=0$$

Hence $a=(1,\,\cos z,\,0).$  By choosing the vector $b$ and $c$ as
follow $$b=(0,1,0),\quad c=(0,0,1)$$ we obtain that $\Upsilon=1.$
Consequently
$$v=\lambda_2(\cos\,z\,\partial_x-\partial_y)-\lambda_3\partial_z$$

The differential equations generated by $v$ and condition \req{80}
in this cases take the form respectively

\begin{equation}
\label{86}
\left\{%
\begin{array}{ll}
\dot{x}=\cos z\,\lambda_2,\\
\dot{y}=-\lambda_2,\\
\dot{z}=-\lambda_3
\end{array}%
\right.   \end{equation}

\begin{equation}\label{87}
  (a,rot{\bold{v}})=0\Leftrightarrow\partial_zp_2-\partial_yp_3+\cos
z\partial_xp_3=0
\end{equation}

After the change $\gamma_1=\sin z\sin x,\quad\gamma_2=\sin z\cos
x,\quad\gamma_3=\cos z$ the system \req{86} and condition \req{87}
can be written as follow
\begin{equation}
\label{88}
\left\{%
\begin{array}{ll}
\dot{\gamma_1}=\frac{I_1}{I_1I_2}\mu_1\gamma_3\\
\\
\dot{\gamma_2}=\frac{I_2}{I_1I_2}\mu_2\gamma_3\\
\\
\dot{\gamma_3}=\frac{-1}{I_1I_2}(I_1\mu_1\gamma_1+
I_2\mu_2\gamma_2)\end{array}%
\right.   \end{equation}

\begin{equation}\label{89}
\sin z(\gamma_3(\frac{\partial \mu_1}{\partial\gamma_2}-
\frac{\partial \mu_2}{\partial \gamma_1})-\gamma_2\frac{\partial
\mu_1}{\partial\gamma_3}+\gamma_1\frac{\partial \mu_2}{\partial
\gamma_3})-\cos x\partial_y\mu_2-\sin
x\partial_y\mu_1=0
\end{equation}
where
\[\left\{%
\begin{array}{ll}
\mu_2=-I_1(\cos x \lambda_3+\sin x \lambda_2),\\
\mu_1=I_2(-\sin x \lambda_3+\cos x \lambda_2)
\end{array}%
\right. \]

 We shall study only the case when

$$ \mu_j=\mu_j(x,z), \,j=1,2$$

 Hence, we obtain the equation \req{84}.

\begin{cor}

{\it The function $\mu_1,\,\mu_2:$

\begin{equation}
\label{90}
\left\{%
\begin{array}{ll}
\mu_1=\displaystyle\frac{\partial
S(\gamma_1,\gamma_2)}{\partial\gamma_1}+\Psi_1(\gamma^2_2+\gamma^2_3,\gamma_1)\\
\mu_2=\displaystyle\frac{\partial
S(\gamma_1,\gamma_2)}{\partial\gamma_2}+
\Psi_2(\gamma^2_1+\gamma^2_3,\gamma_2)
\end{array}%
\right.   \end{equation}

satisfies the equation \req{84}.}
\end{cor}

\begin{cor}

{\it Let $\mu_1,\,\mu_2$ are such that
$$\mu_j=\frac{\partial S(\gamma_1,\gamma_2)}{\partial\gamma_j},\quad j=1,2$$
then the potential function \req{83} and  first integrals \req{85}
are respectively}
\begin{equation}
\label{91}
\left\{%
\begin{array}{ll}
U=\displaystyle\frac{1}{2I_1I_2}(I_1(\displaystyle\frac{\partial
S}{\partial\gamma_1})^2+I_2(\displaystyle\frac{\partial
S}{\partial\gamma_2})^2)-h\\
I_1\omega_1=\displaystyle\frac{\partial S}{\partial\gamma_2},\\
I_2\omega_2=-\displaystyle\frac{\partial S}{\partial\gamma_1},\end{array}%
\right.   \end{equation}
\end{cor}
The following particular cases produces the well known integrable
cases \cite{Koz3}: The Suslov, Kharlamova-Zabelina and Kozlov
subcase.

\textbf{ The Suslov Subcase}

 If
$$S=C_1\gamma_1+C_2\gamma_2,\quad C_j=const,\,j=1,2$$
then
\[\left\{%
\begin{array}{ll}
\mu_1=C_1,\quad \mu_2=C_2\\ U=const. \end{array}%
\right.\] which correspond to the Suslov subcase.

The integration of the equations \req{88} in this case  produces
the following solutions

\[\left\{%
\begin{array}{ll}
\omega_1=\displaystyle\frac{C_2}{I_1},\quad\omega_2=-\displaystyle\frac{C_1}{I_2}\\
\gamma_1=\displaystyle\frac{C_1I_1}{\sqrt{I^2_1C^2_1+I^2_2C^2_2}}\sin\beta\sin
(\displaystyle\frac{\sqrt{I^2_1C^2_1+I^2_2C^2_2}}{I_1I_2} t+\alpha
)+
\displaystyle\frac{I_2C_2\cos\beta}{\sqrt{I^2_1C^2_1+I^2_2C^2_2}}\\
\gamma_2=\displaystyle\frac{C_2I_2}{\sqrt{I^2_1C^2_1+I^2_2C^2_2}}\sin\beta\sin\,(\displaystyle\frac{\sqrt{I^2_1C^2_1+I^2_2C^2_2}}{I_1I_2}
t+\alpha )-\displaystyle\frac{I_1C_1\cos\beta}{\sqrt{I^2_1C^2_1+I^2_2C^2_2}}\\\
\gamma_3=\sin\beta\cos
((\sqrt{\displaystyle\frac{I^2_1C^2_1+I^2_2C^2_2}{I_1I_2}}
t+\alpha)
\end{array}%
\right.\]

 where
$C_1,\,C_2,\,\alpha,\,\beta,\,$ are the arbitrary real constants.

\textbf{ The Kharlamova-Zabelina Subcase }

 If
$$S=\displaystyle\frac{2}{3\sqrt{I_1C^2_1+I_2C^2_2}}(\sqrt{{\tilde{h}+C_1\gamma_1+C_2\gamma_2}})^3+
\displaystyle\frac{CC_2I_2}{{C^2_1I_1+C^2_2I_2}}\gamma_1-
\displaystyle\frac{CC_1I_1}{C^2_1I_1+C^2_2I_2}\gamma_2$$ where
$\tilde{h},\,C_1,\,C_2,\,C$ are arbitrary constants, then

\[\left\{%
\begin{array}{ll}
\mu_1=\displaystyle\frac{C_1}{\sqrt{I_1C^2_1+I_2C^2_2}}\sqrt{\tilde{h}+C_1\gamma_1+C_2\gamma_2}+
\frac{CC_2I_2}{C^2_1I_1+C^2_2I_2}\\
\mu_2=\frac{C_2}{\sqrt{I_1C^2_1+I_2C^2_2}}\sqrt{{\tilde{h}+C_1\gamma_1+C_2\gamma_2}}-\frac{CC_1I_1}{C^2_1I_1+C^2_2I_2}\\
U=\tilde{h}+C_1\gamma_1+C_2\gamma_2
\end{array}%
\right.\]

As a consequence we deduce the Kharlamova-Zabelina subcase
\cite{Kharl}.

The solutions of the equation \req{82}, \req{88} are
\[\left\{%
\begin{array}{ll}
I_1\omega_1=\displaystyle\frac{C_1}{\sqrt{I_1C^2_1+I_2C^2_2}}\sqrt{\tilde{h}+C_1\gamma_1+C_2\gamma_2}+
\displaystyle\frac{CC_2I_2}{C^2_1I_1+C^2_2I_2},\\
I_2\omega_2=-(\displaystyle\frac{C_2}{\sqrt{I_1C^2_1+I_2C^2_2}}\sqrt{{\tilde{h}+C_1\gamma_1+C_2\gamma_2}}-
\displaystyle\frac{CC_1I_1}{C^2_1I_1+C^2_2I_2}),\\
U=\tilde{h}+C_1\gamma_1+C_2\gamma_2\\
\gamma_j=a_j(\tau-C_3)^2+b_j(\tau-C_4)+d_j=\gamma_j (\tau ,C_1,C_2,C_3,C_4),\quad j=1,2\\
\gamma_3=\sqrt{1-\gamma^2_1(\tau ,C_1,C_2,C_3,C_4)-\gamma^2_2
(\tau
,C_1,C_2,C_3,C_4)}\equiv{\sqrt{P_4(\tau,C_1,C_2,C_3,C_4)}}\\
t=t_0+\displaystyle\frac{I_1I_2}{2}\int\displaystyle\frac{d\tau}{\sqrt{P_4(\tau,C_1,C_2,C_3,C_4)}}
\end{array}%
\right.\]

 where
$$a_j=\displaystyle\frac{I_jC_j}{4},\quad
b_j=\displaystyle\frac{CI_1I_2C_1C_2}{C_j(I_1C^2_1+I_2C^2_2)},\quad
d_j=-\displaystyle\frac{\tilde{h}I_jC_j}{I_1C^2_1+I_2C^2_2},$$
 $P_4$ is a polynomial of four
degree in $\tau .$

 \textbf{Kozlov Subcase}

 If we suppose that $I_1=I_2$ and

\[\left\{%
\begin{array}{ll}
S=-2C\arctan{\frac{\gamma_1}{\gamma_2}}+\int D(\gamma^2_1+\gamma^2_2)d(\gamma^2_1+\gamma^2_2)\\
(D(u))^2=\displaystyle\frac{hu^2+\sqrt{1-u}u-C^2}{u^2}\end{array}%
\right.\]

 where $h$ and $C$ are arbitrary real constant, hence,

\[\left\{%
\begin{array}{ll}
\mu_1=-\displaystyle\frac{\gamma_2C}{\gamma^2_1+\gamma^2_2}+\gamma_1D(\gamma^2_1+\gamma^2_2)\\
\mu_2=\displaystyle\frac{\gamma_1C}{\gamma^2_1+\gamma^2_2}+\gamma_2D(\gamma^2_1+\gamma^2_2)\\
U=-h+\sqrt{1-\gamma^2_1-\gamma^2_2}=-h+\gamma_3
\end{array}%
\right.\]

which correspond to the Kozlov subcase.

The equations \req{88} in this case take the form:

\begin{equation}
\label{92}
\left\{%
\begin{array}{ll}
\dot{x}=\displaystyle\frac{C\cos z}{\sin^2z}\\
\dot{y}=\displaystyle\frac{-C}{\sin^2 z}\\
\dot{z}=\displaystyle\frac{(\gamma^2_1+\gamma^2_2)D(\gamma^2_1+\gamma^2_2)}{\sin
z}\end{array}%
\right.   \end{equation}

 which are easy to integrate.

The solutions of the equation of motions are:
\[\left\{%
\begin{array}{ll}
\omega_1=\displaystyle\frac{\gamma_1C}{\gamma^2_1+\gamma^2_2}+\gamma_2D(\gamma^2_1+\gamma^2_2)\\
\omega_2=\displaystyle\frac{\gamma_2C}{\gamma^2_1+\gamma^2_2}-\gamma_1D(\gamma^2_1+\gamma^2_2)\\
x=x_0+{C}\int
\displaystyle\frac{\gamma_3d\gamma_3}{(1-\gamma^2_3)^2D(1-\gamma^2_3)}=x_0+C\int
\displaystyle\frac{\gamma_3d\gamma_3}{\sqrt{(1-\gamma^2_3)P_4(\gamma_3,h,C)}}\\
y=y_0-{C}\int
\displaystyle\frac{d\gamma_3}{(1-\gamma^2_3)^2D(1-\gamma^2_3)}=y_0-C\int
\displaystyle\frac{d\gamma_3}{\sqrt{(1-\gamma^2_3)P_4(\gamma_3,h,C)}}\\
t=t_0+I_1I_2\int\displaystyle\frac{d\gamma_3}{\sqrt{P_4(\gamma_3,h,C)}}\\
P_4(\gamma_3,h,C)\equiv{h\gamma^4_3-2\gamma^3_3-2h\gamma^2_3+2\gamma_3+h-C^2}
\end{array}%
\right.
\]
\begin{cor}
{\it Let $\mu_1,\,\mu_2$ are the functions:
\[\left\{%
\begin{array}{ll}
\mu_1=\Psi_1(\gamma^2_2+\gamma^2_3,\gamma_1)\\
\mu_2=\Psi_2(\gamma^2_1+\gamma^2_3,\gamma_2)\end{array}%
\right.\] then the solutions of \req{88} are the following
functions: }
\end{cor}
\[\left\{%
\begin{array}{ll}
\int\displaystyle\frac{d\gamma_j}{F_j(\gamma_j)}=\frac{I_j}{I_1I_2}(\tau-\tau_0),\quad j=1,2\\
\gamma_3=\sqrt{1-\gamma^2_1(\tau )-\gamma^2_2(\tau )}\\
t=t_0+\int\displaystyle\frac{d\tau}{\sqrt{1-\gamma^2_1(\tau
)-\gamma^2_2(\tau
)}}\end{array}%
\right.\] where
\[\left\{%
\begin{array}{ll}
F_1(\gamma_1)=\Psi_1
(\gamma^2_2+\gamma^2_3,\gamma_1)|_{\gamma^2_2+\gamma^2_3=1-\gamma^2_1}\\
F_2(\gamma_2)=\Psi_2
(\gamma^2_1+\gamma^2_3,\gamma_2)|_{\gamma^2_1+\gamma^2_3=1-\gamma^2_2}\end{array}%
\right.\]

 As a particular case we obtain the
Tisserand Subcase.

\textbf{ Tisserand Subcase}

 The interesting solution of the equation
\req{84} are
\[\left\{%
\begin{array}{ll}
\mu_1=\sqrt{h_1+a_1(\gamma^2_3+\gamma^2_2)+b_1\gamma^2_1+f_1(\gamma_1)}\equiv{\Psi_1( \gamma^2_2+\gamma^2_3,\gamma_1)}\\
\mu_2=\sqrt{h_2+a_2(\gamma^2_3+\gamma^2_1)+b_2\gamma^2_2+f_2(\gamma_2)}\equiv{\Psi_2(
\gamma^2_1+\gamma^2_3,\gamma_2)}\end{array}%
\right.\]
 which produce the following potential function $U:$
$$U=I_1h_1+I_2h_2+(I_1b_1+I_2a_2)\gamma^2_1+(I_1a_1+I_2b_2)\gamma^2_2+(I_1a_1+I_2a_2)\gamma^2_3+I_1f_1(\gamma_1)+
I_2f_2(\gamma_2)$$ where $a_j,\,b_j,\,h_j,\, j=1,2$ are arbitrary
real constants and $f_j,\,j=1,2$ are arbitrary functions.

The case when $f_j(\gamma_j)=\alpha_j\gamma_j,\,j=1,2$ was studied
in \cite{Okuneva}, where $\alpha_j,\,j=1,2$ are real constants.

 The case
when $f_j=0,\,j=1,2$ is well known as Tisserand´s case
\cite{Koz3}.

After integration the equation \req{88} in the Tisserand case we
obtain the following solutions

\[\left\{%
\begin{array}{ll}
I_1\omega_1=\sqrt{h_2+a_2(\gamma^2_3+\gamma^2_1)+b_2\gamma^2_2}\\
I_2\omega_2=-\sqrt{h_1+a_1(\gamma^2_3+\gamma^2_2)+b_1\gamma^2_1}\\
\gamma_1=\sqrt{\frac{h_1+a_1}{a_1-b_1}}\sin
(\sqrt{a_1-b_1}I_1\tau +C_1)=\gamma_1(\tau )\\
\gamma_2=\sqrt{\frac{h_2+a_2}{a_2-b_2}}\sin
(\sqrt{a_2-b_2}I_2\tau +C_2)=\gamma_2(\tau )\\
\gamma_3=\sqrt{1-\gamma^2_1(\tau )-\gamma^2_2(\tau )}\\
t=t_0+I_1I_2\int\frac{d\tau}{\sqrt{1-\gamma^2_1(\tau
)-\gamma^2_2(\tau )}}
\end{array}%
\right.\]

\textbf{  Heavy rigid body in the Veselov case}

In this example we study the problem of non-holonomic dynamics
formulated by Veselov in \cite{Veselov} which in certain sense is
opposite to the Suslov problem. In this problem we consider the
rotational motion of a rigid body around a fixed point and subject
to the non-holonomic constraints
$$(\bold{\gamma},\bold{\omega})\equiv{\dot{y}+\cos
z\dot{x}}=0$$

 Suppose the body rotates in an force field
with potential $U(\gamma_1,\gamma_2,\gamma_3)$. Applying the
method of Lagrange multipliers we write the equations of motion in
the form
\[\left\{%
\begin{array}{ll}
I\dot{\omega}=[I\omega\times{\omega}]+[\gamma\times\displaystyle\frac{\partial
U}{\partial\gamma}]+\lambda\gamma\\
\dot{\gamma}=[{\gamma}\times{{\omega}}] \end{array}%
\right.\]
where $I$ is a matrix such that $I=diag(I_1,I_2,I_3).$

The Descartes approach for this system produces the following
equations:

\[\left\{%
\begin{array}{ll}
\dot{x}=\lambda_2\\
\dot{y}=-\cos z\lambda_2\\
\dot{z}=\lambda_3\end{array}%
\right.\] and

 \[
\displaystyle\frac{\partial p_3}{\partial x}-\frac{\partial
p_1}{\partial z}+\cos z(\frac{\partial p_2}{\partial
z}-\displaystyle\frac{\partial p_3}{\partial y})=0
\] where
\[\left\{%
\begin{array}{ll}
p_1=I_3\sin^2z\lambda_2\\
p_2=(I_3-I_1+(I_1-I_2)\cos^2x)\cos
z\sin^2z\lambda_2+(I_1-I_2)\cos x\sin x\sin  z\lambda_3\\
p_3=(I_2\sin^2x+I_1\cos^2x)\lambda_3+(I_2-I_1)\sin x\cos x\sin
z\cos z\lambda_2\end{array}%
\right.\]

 Finally it is interesting to observe that the construction
the Descartes approach for the Federov case \cite{Fedorov}, i.e.,
$$(\bold{\omega},\bold{\gamma})=a$$
 it is necessary in
 the above example make the change $y\rightarrow{y+at},\,a=const..$
Hence we obtain that that the equations generated by the Descartes
vector field are
\[\left\{%
\begin{array}{ll}
\dot{x}=\lambda_2\\
\dot{y}=-\cos z\lambda_2+a\\
\dot{z}=\lambda_3 \end{array}%
\right.\] and

\begin{equation}\label{93}
  \frac{\partial p_3}{\partial x}-\frac{\partial p_1}{\partial
z}+\cos z(\frac{\partial p_2}{\partial z}-\frac{\partial
p_3}{\partial y}))=0
\end{equation}

 where
\begin{equation}\label{94}\left\{%
\begin{array}{ll}
p_1=I_3\sin^2z\lambda_2\\
p_2=(I_3-I_1+(I_1-I_2)\cos^2x)\cos
z\sin^2z\lambda_2+(I_1-I_2)\cos x\sin x\sin  z\lambda_3+\\
a((I_1\sin^2 x+I_2\cos^2 x)\sin^2{z}+I_3\cos^2{z})\\
p_3=(I_2\sin^2x+I_1\cos^2x)\lambda_3+(I_2-I_1)\sin x\cos x\sin
z\cos z\lambda_2
\end{array}%
\right.
\end{equation}

 The equation \req{93} can be represented as follow
\begin{equation}\label{95}\left\{%
\begin{array}{ll}
sin^2 z(I_3\sin^2 z+\cos^2 z (I_1\sin^2 x+I_2\cos^2x)\partial_z\lambda_2(x,z)+ \\
(I_2\sin^2x+I_1\cos^2x)\partial_x\lambda_3(x,z)+\\
\cos x\sin\,x\sin z\cos\,
z(I_1-I_2)(\partial_z\lambda_3(x,z)-\partial_x\lambda_2(x,z))+\\
 \sin z\cos\,z(3(I_1-I_2)\sin^2z\cos^2x+
3I_1-I_3)\sin^2z+I_1+I_2)\lambda_2(x,z)\\
\cos x\sin\,x(-2+\cos^2z)(I_1-I_2)\lambda_3(x,z)+ 2a\sin
z\cos^2z(I_1\sin^2x+I_2\cos^2x)=0.
\end{array}%
\right.
\end{equation}
\begin{prop}

{\it Let $\lambda_2,\,\lambda_3$ are the solutions of the linear
partial differential equations \req{95} then the functions}
\begin{equation}\label{96}
\left\{%
\begin{array}{ll}
\omega_1=\gamma_2\frac{\lambda_3}{\sin
z}-\gamma_1\gamma_3\lambda_2-a\gamma_1\\
\omega_2=-\gamma_1\frac{\lambda_3}{\sin
z}-\gamma_2\gamma_3\lambda_2-a\gamma_2\\
\omega_3=\sin^2z\lambda_2-a\gamma_3
\end{array}%
\right.\end{equation}{\it are the first integral of the equations
of the rigid body in the Veselov-Fedorov case.} \end{prop}

In particular if $$I_1=I_2,$$ then the solutions of the above
equation are \quad
\[\left\{%
\begin{array}{ll}
\lambda_3=\sin z\,C(z)\\
 \sqrt{I_3\sin^2 z+\cos^2 z I_2}\sin^2
z\lambda=a\Omega (z)+K(x)\\
 \end{array}%
\right.
\]where $C(z)$ and $K(x)$ are arbitrary
functions and
$$\Omega (z)\equiv{\int\displaystyle\frac{I_2\sqrt{(I_3\sin^2
z+\cos^2 z
I_2})(2\sin\,z-\sin\,5z+\sin\,3z)dz}{(-I_3\cos\,4z+4I_3\cos\,2z-3I_3-I_2+I_2\cos\,4z)}
}$$
 Hence, from \req{96} we obtain

\begin{equation}\label{97}
\left\{%
\begin{array}{ll}
\omega_1=\gamma_2C(z)-\gamma_1\gamma_3\displaystyle\frac{a\Omega
(z)+K(x)}{\sqrt{I_3\sin^2 z+\cos^2 z I_2}\sin^2
z}-a\gamma_1\\
\omega_2=-\gamma_1C(z)-\gamma_2\gamma_3\displaystyle\frac{a\Omega
(z)+K(x)}{\sqrt{I_3\sin^2 z+\cos^2 z I_2}\sin^2
z}-a\gamma_2\\
\omega_3=\displaystyle\frac{a\Omega (z)+K(x)}{\sqrt{I_3\sin^2
z+\cos^2 z I_2}}-a\gamma_3
\end{array}%
\right.\end{equation} Thus,  we easily deduce the relation
$$(I_3\sin^2 z+\cos^2 z I_2)(\omega_3+a\gamma_3)^2=(a\Omega (z)+K)^2.$$
In particular if $a=0$ and $K(x)=C_1=const$ then  we obtain   the
well known first integral in the Veselov case \cite{Borisov}

\subsection*{Acknowledgments}
This work was partly supported by the Spanish Ministry of
Education through projects DPI2007-66556-C03-03,
TSI2007-65406-C03-01 "E-AEGIS" and Consolider CSD2007-00004
"ARES".



\begin{thebibliography}{99}
\bibitem{Bates}
Bates,L. and Sniatycki J.,  Nonholonomic reduction. {\em Rep.
Math. Phys.} \textbf{32} (1), 1992, 99-115
\bibitem{Koz1}
Kozlov, V.V.,  Dynamical system X, General theory of
 vortices. Spriger, ( 2003).

\bibitem{Ram1}
 Ramirez, R. and Sadovskaia N., Descartes
 approach for Non-holonomic systems,
  {\em J. Math. Sciences (N.Y.)}, vol. 128, Number 2
 (2005), 2812-2817.
\bibitem{Borisov}
Borisov A.V. and Mamaev I.S., Strange attractors in rattleback
dynamics.,  {\em  Phys. Usp.} 46 (4), (2003), pp.393-403.
\bibitem{Dub} Duboshin, G.H., Celest. Mech., Ed.  Nauka,1968,
Moscow  (in Russian).
\bibitem{Kharl}
E.I. Kharlamov-Zabelina,  Rapid rotation of a rigid body about a
fixed point under the presence of a non-holonomic constraints.
{\em Vestnik Moskovsk. Univ., Ser. Math. Mekh., Astron.,
Fiz.Khim.}, V.6, (1957), pp.25-34 (in Russian).

\bibitem{Koz2}
Kozlov V V, Several problems on dynamical systems and mechanics,
{\em Nonlinearity } {\bf 21} (2008) T149-T155.

\bibitem{Ros}
 Rosenberg, R.,
Analytical Dynamics, Ed. Heinemann, New York, 1977.
\bibitem{Ram2} Ram\'{\i}rez Rafael and Sadovskaia N., On the
 dynamics of non-holonomic systems, {\em Reports on Mathematical
 Physics,} vol. 60 (2007).
\bibitem{Ram3}Ram\'{\i}rez Rafael and Sadovskaia N., Inverse Problem in
Celestial Mechanic,{\em Atti. Sem. Mat.Fis. Univ. Modena e Reggio
Emilia,} LII, (2004) pp.47-68 .
\bibitem{Okuneva}
Okuneva, G.G.,  Dvijenie tverdogo tela s nepodvijnoy tochkoi pod
deistviem negolonomnoi sviazi v niutonovskom pole. {\em Mekhanika
tverdogo tela,} Kiev, {\bf 18}, (1986),  pp.40-43 (in Russian).

\bibitem{Koz3}
Kozlov, V.V., Symmetries, Topology and Resonances in Hamiltonian
Mechanics. Springer-Verlag. Berlin,
\bibitem{Veselov}
Veselova, L.E., Novie sluchai integriruemosti uravneniy dvijenia
tverdogo tela pri nalichii negolonomnoi sviazi. {\em  Sbornik
Geometria, dif. uravnenia i mekhanika, MGU,} (1986), pp.64-68 (in
Russian).

\bibitem{Fedorov}
Fedorov I N, O dvuh integriruemih negolonomnih sistemah v
klasicheskoi mehanike,{\em Bectnik MGU, serie Math. Mech.} {\bf
4}, (1989), pp.38-41.
\bibitem{NF} Neimark Ju.I. and Fufaev N.A., Dynamics of Nonholonomic Systems, American
Mathematical Society, (1972).

\bibitem{Sad} Sadovskaia, N.,  Inverse problem in theory of ordinary
differential equations. Thesis Ph. Dr., Univ.  Polit\'ecnica  de
Catalu\~na , (2002) (in Spanish)





\end{thebibliography}
\end{document}